 {\everymath{\displaystyle\everymath{}}\array}%
 {\endarray}

\documentclass[10pt,a4paper]{article}
\usepackage{amssymb}
\newtheorem{thm}{Theorem}
\newtheorem{defn}{Definition}
\newtheorem{cor}{Corollary}
\newtheorem{lem}{Lemma}

\newtheorem{prop}{Proposition} 
\newtheorem{ex}{Example}
\newtheorem{rk}{Remark}

\usepackage{amssymb}
\usepackage{color}
\usepackage{colordvi}
\usepackage{amsmath}

\newcommand{\field}[1]{{\mathbb #1}}
\newcommand{\R}{\field{R}}
\newcommand{\C}{\field{C}}

\newcommand{\Z}{\field{Z}}

\newcommand{\N}{\field{N}}

\newcommand{\T}{\field{T}}
\newcommand\cutoffsum{\mathop{-\hskip -4mm\sum}\limits}

\def \endsquare{ $\sqcup \!\!\!\! \sqcap$ }
\def\up{{\underline p}}
\def\un{{\underline n}}
\def\ux{{\underline x}}
\def\uy{{\underline y}}

\def\cutoffint{-\hskip -10pt\int}

\def\otherterm#1{{\it#1}}
\def \e {{\epsilon}}

\def \Ci {{C^\infty}}
\def \l {{\lambda}}

\def \Cl {{C\ell}}

\begin{document}

\title{\bf Discrete sums of classical symbols on $\Z^d$ and zeta functions
  associated with Laplacians on tori}

\author{  Sylvie
PAYCHA }
\maketitle 

\section*{Abstract\footnote{AMS classification: 11E45, 58C40, 47G30}}{\it We
  prove the uniqueness of a  
  translation invariant extension
  to non integer order classical symbols of the ordinary discrete sum on
  $L^1$-symbols, which  we then describe
  using an Hadamard finite part procedure for sums over 
  integer points of infinite unions of nested convex polytopes in $\R^d$.
   This canonical regularised sum is the building block to construct
   meromorphic extensions of the ordinary sum on holomorphic symbols.  Explicit formulae for
  the complex residues at their  
  poles are given  in terms of noncommutative residues of classical symbols, thus 
  extending results of  Guillemin,  Sternberg and Weitsman. These formulae
  are then applied to zeta functions associated with quadratic forms and with  Laplacians on tori. }
\section*{Introduction}
Just as the ordinary $L^2$-trace on trace-class classical pseudodifferential
operators on a closed manifold is known not to extend to a trace to the whole
algebra of classical operators, one does not expect the   Riemann integral
(resp.discrete sum)  on $L^1$ classical pseudodifferential symbols on $\R^d$ to  extend to an $\R^d$ (resp. $\Z^d$)
-translation invariant linear form on the whole algebra of symbols. \\
We show that these nevertheless  have a canonical  $\R^d$ (resp. $\Z^d$)
-translation invariant linear extension\footnote{By linear we mean that it
  preserves linear combinations that lie in the set.} to the set of {\it non integer order}
classical symbols  on $\R^d$ (see Theorem \ref{thm:uniqueness}), namely the
canonical regularised integral $\cutoffint_{\R^d}$  (resp. the canonical
regularised discrete sum $\cutoffsum_{\Z^d}$)
\footnote{We borrow the 
notation $\cutoffint_{\R^d}$  from \cite{L} and transpose it to the discrete
sum but warn the reader that the same notation is used for the Dixmier trace
by A.Connes.}. \\
It is known (see e.g. \cite{L}) that the
canonical regularised integral $\cutoffint_{\R^d}$ can be expressed  as an
Hadamard finite part of ordinary integrals over euclidean balls of radius tending to
infinity and that its translation invariance is related to a Stokes property
\cite{P2}. Integrating this canonical regularised integral over a closed
manifold gives rise to the unique extension of the $L^2$-trace defined on trace-class classical
pseudodifferential operators to the set of non integer order classical
pseudodifferential operators \cite{MSS}, namely to the   canonical trace introduced by Kontsevich and Vishik \cite{KV}.
\\ Less is known about  extensions of the ordinary discrete sum.
Whereas one does not expect  the canonical
regularised discrete sum to vanish on derivatives, its value on derivatives is entirely determined by
its restriction to $L^1$ symbols (see Proposition \ref{prop:extrho}). We express the  canonical regularised  sum as an Hadamard finite
part for expanded convex polytopes with increasing size, using  a
generalisation of the  Euler-MacLaurin formula in order to compare these with the corresponding integrals
over these polytopes.\\
  The Khovanskii-Pukhlikov formula \cite{KP} generalises the Euler-MacLaurin
  formula   to  polynomials  in higher dimensions in as far as it compares
the discrete sum  over integer points  with the
integral  of a polynomial function on  expanded convex polytopes which are
integral and regular\footnote{It was then generalised to simple integral 
polytopes by Cappell and Shaneson \cite{CS} and subsequently by Guillemin
 \cite{G2} and Brion and
Vergne \cite{BV}; these generalisations   involve corrections
to the  Khovanskii-Pukhlikov formula when the simple polytope is not regular.}. Guillemin,  Sternberg and Weitsman \cite{GSW}   extended  the
Khovanskii-Pukhlikov formula to classical symbols.  Whereas it
is an exact formula for polynomials, it becomes an asymptotic formula in the
case of classical symbols and a new constant term $C(\sigma)$ depending on the
symbol $\sigma$ 
arises,  which vanishes for
polynomials. For gauged symbols $z\mapsto \sigma(z)$, the authors showed that
$z\mapsto C(\sigma(z))$  is
holomorphic and interpreted  $C(\sigma)$ as a ``regularised
version'' of the difference between the infinite sum  ''$\sum_{\un\in
  \Z^d}\sigma(\un)$'' of the symbol  over integer points of $\R^d$ and the
infinite integral ``$\int_{\R^d} \sigma(\ux)\, d\ux$''. \\
In this paper, using the canonical regularise sum $\cutoffsum_{\Z^d}$, we refine and generalise some of their results to any local holomorphic
perturbation $\sigma(z)$  of $\sigma=\sigma(0)$ with non constant affine order
$\alpha(z)$  showing
the following statements.
\begin{enumerate}
\item The map $z \mapsto \sum_{\un\in \Z^d}\sigma(z)(\un)$ which is holomorphic on a half 
plane Re$(\alpha(z))<-d$, 
  extends to a meromorphic map $z\mapsto \cutoffsum_{\un\in \Z^d}\sigma(z)(\un)$
  on the whole complex plane with simple poles.  The pole at $z=0$  is proportional to
  the noncommutative residue of $\sigma$ (see Theorem  \ref{thm:maintheorem})
$${\rm Res}_{z=0}\cutoffsum_{\un\in
  \Z^d}\sigma(z)(\un)=-\frac{1}{\alpha^\prime(0)} {\rm res}(\sigma).$$
\item When $\sigma$ has non integer order, the constant term in the Laurent
  expansion $ \sum_{\un\in \Z^d}\sigma(z)(\un)$ around $z=0$ coincides with
  the canonical regularised sum ${\rm fp}_{z=0} \cutoffsum_{\un\in \Z^d}\sigma(z)(\un)=
  \cutoffsum_{\un\in \Z^d}\sigma(\un)$. 
\item In general, the finite parts ${\rm fp}_{z=0}\cutoffsum_{\un\in
    \Z^d}\sigma(z)(\un)$ for different holomorphic perturbations $\sigma\mapsto
  \sigma(z)$ only differ if
  the perturbations do not coincide on the unit sphere, in which case they
  differ  by a term
  proportional to some noncommutative residue   (see Corollary \ref{cor:maincorollary}).
\end{enumerate}
These results strongly rely on corresponding known results \cite{KV} for
 regularised integrals of symbols which we recall in an
Appendix, namely the fact that the ordinary integral
$z \mapsto \int_{\R^d}\sigma(z)(\ux)\, d\ux$ has a meromorphic extension $z
\mapsto \cutoffint_{\R^d}\sigma(z)(\ux)\, d\ux$ to the
whole complex  plane with simple poles and that the pole at $z=0$ is
proportional
to the noncommutative residue of $\sigma$ (see Theorem \ref{thm:KVsymbol} in
the Appendix). \\
In the last part of the paper (sections 6 and 7), we apply our results to the zeta
function $$Z_q(s):={\rm fp}_{z=0}\cutoffsum_{\Z^d-\{0\}} q( \un)^{-s}\,
\vert \un\vert^{-z}$$
associated with a
positive definite quadratic form $q$ on $\R^d$  (see  Theorem 
\ref{thm:quadraticzeta}), resp. to the zeta function $$  \zeta_{\Delta_d}(s):={\rm
  fp}_{z=0}\cutoffsum_{\Z^d-\{0\}} \vert \un\vert^{-2s+z}$$ associated with the Laplacian $\Delta_d$ on the
$d$-dimensional torus (see Proposition \ref{prop:zetaDelta}). Here $\vert \ux\vert$ stands for the euclidean norm on $\R^d$.\\ At arguments  $s$ with non positive real
part, $Z_q(s)$, resp. $\zeta_{\Delta_d}(s)$ coincides with  the constant
$C(q^{-s})$ (see Proposition \ref{prop:zetaqhol}), resp. $C(\ux\mapsto \vert \underline
x\vert^{-2s})$, leading  to a formula for
 $\zeta$-determinants of Laplacians on tori (see Proposition
 \ref{prop:zetadeterminant})
$${\rm det}_\zeta (\Delta_d)= \exp\left( -{\partial_s}_{\vert_{ s=0^-}}C(\ux\mapsto \vert
  \ux\vert^{-2s})\right)$$
where the $ {\partial_s}_{\vert_{ s=0^-}}$  stands for the derivative 
in the half plane Re$(s)\leq 0$. \\ \\
The paper is organised as follows.
\begin{enumerate}
\item Sums of classical symbols on positive integers
\item  A canonical regularised integral and discrete sum on non integer order symbols
\item Concrete realisations of the canonical regularised integral and discrete sum
\item Sums of holomorphic symbols on $\Z^d$
\item Zeta functions associated with quadratic forms
\item Zeta functions associated with Laplacians on tori
\item Appendix: Prerequisites on regularised integrals of symbols
\end{enumerate}
\section*{Acknowledgements}
I would like to thank Victor Guillemin, Shlomo Sternberg and  Jonathan 
Weitsman for their 
comments on a preliminary version of this paper. Let me also 
thank Viatcheslav Kharlamov for his comments following a talk at the Max
Planck institute in Bonn during which I 
reported on
these results   as well as the institute itself,  where this article was completed. 
\vfill \eject \noindent
\section{Discrete sums of classical symbols on integers} 
 Before investigating discrete sums on polytopes, let us first consider discrete sums of classical symbols
on ordinary  integers. Along the lines of \cite{MP}  we use the Euler-MacLaurin formula which relates a sum to an
integral. It involves the
 Bernouilli numbers  defined by the generating series:
\begin{equation} \label{eq:Bernouillinumb}
\sum_{n=0}^\infty B_n \frac{t^n}{n!} = \frac{t}{e^t-1}.
\end{equation} 
Bernouilli polynomials are defined similarly:
 \begin{equation} \label{eq:Bernouillipoly}
\sum_{n=0}^\infty B_n(x) \frac{t^n}{n!} = \frac{t\, e^{t\, x}}{e^t-1},
\end{equation}
so that in particular,  $ B_n(0)=B_n$. For example
 $B_1(x)= -\frac{1}{2}+ x.$ \\  The  Euler-MacLaurin formula
 states that for any   $ f\in \Ci(\R)$  and any two integers $M<N$  (see e.g.  \cite{Ha})
\begin{eqnarray}\label{eq:EulerMacLaurin} \sum_{n=M}^{N} f(n)&=&\frac{f(M)+f(N)}{2} +\int_1^N f(x)\, dx+
\sum_{k=2}^K (-1)^{k} \frac{B_k}{k!}\left( f^{(k-1)}(N)-f^{(k-1)}(M)\right)\nonumber\\
&+&
\frac{(-1)^{K-1}}{K!} \int_M^N \overline{B_{K}} (x)\, f^{(K)}(x) \, dx
\end{eqnarray} 
with $\overline{B_k}(x)= B_k\left(x-[x] \right)$ and where $K$ is any positive
integer larger than $1$. 
 \\ \\
When $f$ is polynomial of degree $d$, this reduces to 
$$  \sum_{n=M}^{N} f(n)=\frac{f(M)+f(N)}{2} +\int_M^N f(x)\, dx+
\sum_{k=2}^{d+1} (-1)^{k} \frac{B_k}{k!}\left( f^{(k-1)}(N)-f^{(k-1)}(M)\right). 
$$ 
 As it was shown in  \cite{GSW}, see also \cite{MP},   this formula
 gives interesting information when applied to a  symbol $\sigma$ in the
 algebra  $
 CS_{\rm c.c}(\R)$ of classical pseudodifferential symbols with constant coefficients defined in formula (\ref{eq:CSRn}) of the Appendix.  In particular,  the Euler-MacLaurin formula
  provides a  control on the asymptotics
  of  $\sum_{n=M}^{N} \sigma(n)$ as $M=-N\to -\infty $ with $N\to \infty$
  \begin{eqnarray*}  \sum_{n=-N}^{N} \sigma(n)&\sim_{N\to \infty}&\frac{
  \sigma(-N)+f(N)}{2} +\int_{-N}^N f(x)\, dx+
\sum_{k=2}^K (-1)^{k} \frac{B_k}{k!}\left( \sigma^{(k-1)}(N)-\sigma^{(k-1)}(-N)\right)\nonumber\\
&+&
\frac{(-1)^{K-1}}{K!} \int_\R \overline{B_{K}} (x)\, \sigma^{(K)}(x)
 \, dx.
\end{eqnarray*} Picking the constant term we
   define an Hadamard finite part
   $ {\rm fp}_{N\mapsto \infty} \sum_{n=-N}^{N}
  \sigma (n)$ (see also \cite{MP}) which relates to  the finite part ${\rm fp}_{N\to \infty} \int_{-N}^N
  \sigma(x)\, dx$ using (\ref{eq:EulerMacLaurin}): 
  \begin{eqnarray}\label{eq:EulerMacLaurinreg}  {\rm fp}_{N\mapsto \infty} \sum_{n=-N}^{N} 
  \sigma(n)&=& 
{\rm fp}_{N\mapsto \infty}  \cutoffint_{-N}^N \sigma(x)\, dx+ \frac{{\rm fp}_{N\to \infty}\sigma(-N)+{\rm fp}_{N\to \infty} \sigma(N)}{2} \nonumber\\
  &+&
\sum_{k=2}^K (-1)^{k} \frac{B_k}{k!}\left( {\rm fp}_{N\to \infty}
 \sigma^{(k-1)}(N)-{\rm fp}_{N\to \infty}\sigma^{(k-1)}(-N)\right)\nonumber\\
&+&
\frac{(-1)^{K-1}}{K!} \int_\R\overline{B_{K}} (x)\,
 \sigma^{(K)}(x) \, dx
\end{eqnarray} 
where this last term is actually a convergent integral. \\
The difference
 ${\rm fp}_{N\mapsto \infty} \sum_{n=-N}^{N} 
  \sigma(n)- 
{\rm fp}_{N\mapsto \infty}  \cutoffint_{-N}^N \sigma(x)\, dx$ involves finite parts ${\rm fp}_{N\to \infty}
 \sigma^{(j)}(N)$ and ${\rm fp}_{N\to \infty}
 \sigma^{(j)}(-N)$ which  vanishes  when $\sigma$ has non integer order. On
 non integer order

 of the  Appendix.  symbols,  the Hadamard finite parts ${\rm fp}_{N\mapsto
   \infty}  \cutoffint_{-N}^N \sigma(x)\, dx$, resp.  ${\rm fp}_{N\mapsto \infty} \sum_{n=-N}^{N} 
  \sigma(n)$  turn out to be independent of the choice of
 parameter $N$ and  give rise to $\R$ (resp. $\Z$) -translation invariant extensions
 $\sigma\mapsto \cutoffint_{\R} \sigma(x)\, dx:={\rm fp}_{N\mapsto
   \infty}  \cutoffint_{-N}^N \sigma(x)\, dx$, resp. $ \sigma\mapsto\cutoffsum_{\Z} \sigma(n):={\rm fp}_{N\mapsto \infty} \sum_{n=-N}^{N} 
  \sigma(n) $ of the
 ordinary integral (resp. discrete sum) on $L^1$ symbols to the set $CS_{\rm
   c.c}^{\notin \Z}(\R)$ of non integer order symbols on $\R$ defined in
 (\ref{eq:CSnoninteger})  of the Appendix:
$$ \sigma\in
CS^{\notin \Z}_{\rm c.c}(\R) \Longrightarrow\left( \cutoffint_{\R} t_p^*\sigma(x)\, dx= \cutoffint_{\R} \sigma(x)\, dx\quad
\forall p\in \R; \quad
\cutoffsum_{\Z} \sigma(n)=\cutoffsum_{\Z} t_p^*\sigma(n)\quad \forall p\in
\Z\right)$$  where
  $t_p^*\sigma:= \sigma(p+\cdot)$. These are related via the map $\sigma\mapsto
 C(\sigma)$:
   $$ C (\sigma):=\cutoffsum_{n\in \Z}
  \sigma(n)- 
  \cutoffint_\R \sigma(x)\, dx
  =
\frac{(-1)^{K-1}}{K!} \int_\R \overline{B_{K}} (x)\,
 \sigma^{(K)}(x) \, dx$$
which is clearly  translation invariant on non integer order symbols  i.e.
$C\left( t_p^*\sigma\right)= C(\sigma)$ for any $p\in \Z$. \\
  We apply  this to the symbol $\sigma_s(x)= \chi(x)\, x^{-s}$ where $\chi$ is a
  smooth cut-off function which vanishes in a small neighborhood of $0$ and is
  identically one outside the unit interval. This symbol is of non integer order when
  $s\notin \Z$ and we have
  $$\cutoffsum_{n\in \Z-\{0\}}
  n^{-s}= 
  \cutoffint_\R \sigma_s(x)\, dx+C (\sigma_s)\quad \forall s\notin \Z.$$
  The map  $s\mapsto C(\sigma_s)$ turns out to  be holomorphic so that  
  $\lim_{s\to s_0} C (\sigma_s)= C(\sigma_{s_0})$ and the Riemann zeta function
  reads:
  $$\zeta(s_0):= {\rm fp}_{s=s_0}\cutoffsum_{n=1}^{\infty} 
  n^{-s}= \frac{1}{2}{\rm fp}_{s=s_0}\cutoffsum_{n\in \Z-\{0\}} 
  \sigma_s(n)=\frac{1}{2}\,
 {\rm fp}_{s=s_0} \cutoffint_\R \sigma_{s} (x)\, dx+\frac{1}{2}\, C (\sigma_{s_0})\quad \forall s_0\in
  \C.$$
For Re$(s)\leq 0$, the expression $\overline\sigma_s(x)= x^{-s}$ is defined on non
non negative real numbers and  for any complex number $s_0$ with non negative real
part, we have $\zeta(s_0)=\frac{1}{2}\, C(\overline\sigma_{s_0})$. 
 
  \begin{rk}Note that at integer points $s_0$, one is to expect that
  $${\rm fp}_{s=s_0}\cutoffsum_{n=1}^{\infty} 
  n^{-s}\neq \cutoffsum_{n=1}^{\infty} 
  n^{-s_0}.$$ In particular, taking $s_0=-(2k+1)$ with $k$  a non negative integer we have
  $$ \zeta(-(2k+1))= -\frac{B_{2k+2}}{2k+2}\quad \neq\quad  \cutoffsum_{n=1}^{\infty} 
  n^{2k+1}=0.$$
  \end{rk}
The   Laplacian $\Delta_{S^1}$ on the unit circle
 $S^1\simeq \R/2\pi \Z$  induced by the Laplacian $ -\frac{d^2}{dt} $ on $\R$ has purely discrete
 spectrum given by $\{\l_n= n^2, n\in \Z\}$ and for Re$(s)>1$ we have 
 $$\zeta(s)=\frac{1}{2}\,  \sum_{  n\in \Z-\{ 0\}}\lambda_n^{\frac{s}{2}}=
 \frac{1}{2}\, {\rm tr}\left(\Delta_{S^1}^{-\frac{s}{2}}\right).$$
   Laplacians $\Delta_{\T^d}$ on $d$-dimensional tori $\T^d\simeq \R^d/ 2\pi
   \Z^d$ induced by the Laplacian $\Delta_d:= -\sum_{i=1}^d \frac{d^2}{dx_i}$ on $\R^d$ 
   have purely discrete
 spectrum given by $\{\l_{\un}= \vert \un\vert^2, \un\in \Z^d\}$  where $\vert
 \cdot \vert$ is the euclidean norm on $\R^d$. Zeta functions associated with
 Laplacians on tori  yield a natural generalisation of the Riemann zeta
 function:
 $$\zeta_{\Delta_{d}}(s)=
 \sum_{\un\in \Z^d-\{ 0\}}\lambda_{\un}^{-s}=
 {\rm tr}\left(\Delta_{\T^d}^{-s}\right)\quad {\rm for} \quad {\rm Re}
 (s)>\frac{d}{2}.$$
 The subsequent constructions lead to a   generalisation of    the one dimensional results described above to
 higher dimensions providing along the way a meromorphic extension $z\mapsto
 \cutoffsum_{\un\in \Z^d-\{ 0\}}\lambda_{\un}^{-s-z/2}$ of $\sum_{\un\in
   \Z^d-\{ 0\}}\lambda_{\un}^{-s}$ defined for ${\rm Re}(s)>>0$.
 To carry out this generalisation we strongly rely on the  symbolic feature of the eigenvalues $\lambda_{\un}$ of $\Delta_{d}$ which is rather specific to these operators; a 
 natural question   would be to generalise some of these results to any
  elliptic pseudodifferential operator on a closed manifold with positive
   order and positive leading symbol, which is  probably a difficult issue to
   handle.
\vfill \eject \noindent
\section{A canonical regularised  integral and discrete sum on non integer
  order symbols}
\subsection{Translation invariant linear forms on non integer order symbols}
 Given a symbol $\sigma$ in the algebra $ CS_{\rm c.c}(\R^d)$ of classical
 pseudodifferential symbols with constant coefficients defined in the
 Appendix, and a vector $\up\in \R^n$,
  the map $\up \mapsto t_{ \up}^* \sigma:= \sigma(\cdot+  \up)$ has the following Taylor expansion
  at $\up=0$: 
\begin{equation}\label{eq:Taylorexpansion}t_{\up}^* \sigma(\ux):= \sum_{\vert \beta\vert\leq N}\partial^\beta \sigma (\ux)\,
\frac{\up^\beta}{\beta!} +\sum_{\vert \beta\vert=N+1}\frac{
\up^\beta }{\beta!} \,
 \, \int_0^1 (1-t)^N\, \partial^\beta  \sigma (\ux+ t\up)\,  dt\quad
\forall \ux\in \R^d.
\end{equation}
The Riemann integral and discrete sum induce two linear forms  $\int_{\R^d}: CS_{\rm c.c}^{<-d}(\R^d)\to \C$ and 
$\sum_{\Z^d}: CS_{\rm c.c}^{<-d}(\R^d)\to \C$ on the subalgebra $ CS_{\rm
  c.c}^{<-d}(\R^d)$ of classical pseudidifferential symbols the order of which
has real part $<-d$ (defined in the Appendix) which satisfy the following translation
invariance property:
$$\int_{\R^d} t_{\up}^*\sigma= \int_{\R^d} \sigma\quad \forall \up\in \R^d;
\quad \sum_{\Z^d} t_{\up}^*\sigma= \sum_{\Z^d} \sigma\quad \forall \up\in
\Z^d$$ so that
$$\sum_{0<\vert \beta\vert\leq N} \left(
\int_{\R^d} \partial^\beta \sigma (\ux)\,d\ux\right)
\frac{\up^\beta}{\beta!} +\sum_{\vert \beta\vert=N+1}\frac{
\up^\beta}{\beta!} \, \, \int_0^1 (1-t)^N\, \left( \int_{\R^d} \partial^\beta\sigma (\ux+
  t\up)\,d\ux\right)\,   dt=0\quad\forall \up\in \R^d,$$ and similarly
replacing the integral by a discrete sum over $\Z^d$ and letting $\up$ vary in
$\Z^d$. 
Differentiating  this identity in the components of
$\up$ at $\up=0$  we infer  that $$\int_{\R^d} \partial^\beta \sigma
(\ux)\,d\ux=0\quad \forall \beta \neq 0,$$
and thereby recover the fact that the ordinary integral vanishes on derivatives of symbols of order
$<-d$. \\ \\
Let us further extend  the notion of translation invariant
linear form\footnote{By linear,as observed in the introduction,  we mean that it
  preserves linear combinations that lie in the set.} to the set $CS^{\notin \Z}(\R^d)$ of non integer order symbols defined in formula
(\ref{eq:CSnoninteger}) of the Appendix. To do so, we observe that for $\vert \beta\vert $ large enough, 
 the map $\ux\mapsto \partial^\beta\sigma (\ux+
  t\up)$ lies in $L^1(\R^d)$.
\begin{defn}Given a linear form $\rho: L^1(\R^d) \to \C$, such that 
\begin{equation}\label{eq:rho}\rho\left(t_{\up}^* \sigma\right)= \sum_{\vert \beta\vert\leq N}\rho(\partial^\beta \sigma )\,
\frac{\up^\beta}{\beta!} +\sum_{\vert \beta\vert=N+1}\frac{\up^\beta }{\beta!} \,
\, \int_0^1 (1-t)^N\,\rho\left( \partial^\beta  \sigma (\cdot+ t\up)\right)\,  dt\quad
\forall \up\in \R^d, \end{equation}
  a linear form $\lambda: CS_{\rm c.c}^{\notin \Z}(\R^d)\to \C$ which
  restricts to $\rho $ on  $CS_{\rm c.c}^{\notin \Z}(\R^d)\cap CS_{\rm c.c}^{<-n}(\R^d) $
\begin{itemize}
\item is  extended
  to translated symbols by:
\begin{eqnarray*}
t_{\up}^*\lambda\left( \sigma\right)&:=& \lambda\left( t_{\up}^* \sigma\right)\\
&=&\sum_{\vert \beta\vert\leq N} \lambda \left(
 \partial^\beta \sigma \right)
\frac{\up^\beta}{\beta!} +\sum_{\vert \beta\vert=N+1}\frac{\up^\beta }{\beta!} \,
\, \int_0^1 (1-t)^N\, \rho \left( \partial^\beta\sigma (\cdot+
  t\up)\right)\,   dt,
\end{eqnarray*}
for any large enough  integer $N$. 
\item  The linear form $\lambda$ is said
  to be $\R^d$ (resp. $\Z^d$) -translation invariant whenever for any large
  enough integer $N$, for any $ \up\in \R^d$  (resp.  for any $
\up \in \Z^d$) we have:
$$\sum_{0<\vert \beta\vert\leq N}
\lambda( \partial^\beta \sigma) 
\frac{\up^\beta}{\beta!} +\sum_{\vert \beta\vert=N+1}\frac{\up^\beta }{\beta!} \,
 \, \int_0^1 (1-t)^N\, \rho \left(  \partial^\beta\sigma (\cdot+
  t\up)\right)\,   dt=0.$$
\end{itemize}
\end{defn}
\begin{rk}  These definitions  are   independent of the choice of the integer $N$
  as a result of (\ref{eq:rho}) and the fact that $\lambda$ coincides with
  $\rho$ on $CS_{\rm c.c}^{<-d}(\R^d)$. 
\end{rk}
The following proposition shows that an $\R^d$ or $\Z^d$-translation
invariant linear extension of a linear form $\rho$  on $L^1$-symbols to non integer
order symbols, is entirely determined by $\rho$. 
\begin{prop} \label{prop:extrho} Let $\lambda: CS_{\rm c.c}^{\notin \Z}(\R^d)\to \C$ be a
  linear form   which restricts  on 
  $CS_{\rm c.c}^{\notin \Z}(\R^d)\cap CS^{<-n}(\R^d) $  to a
  linear form $\rho: L^1(\R^d)\to \C $ satisfying (\ref{eq:rho}).
 \begin{enumerate}
\item  If $\lambda$ is  $\R^d$-translation invariant, it  vanishes on
  partial derivatives. 

\item If $\lambda$ is  $\Z^d$-translation invariant, its value on derivatives
  $\lambda(\partial^\beta \sigma),\quad  \beta \neq 0, \sigma\in CS_{\rm c.c}^{\notin \Z}(\R^d)$  is entirely
  determined by the restriction $\rho$ to $CS^{<-n}(\R^d) $.
\end{enumerate}  
\end{prop}
{\bf Proof:}\begin{enumerate}
\item  Differentiating  the identity$$\sum_{0<\vert \beta\vert\leq N}
\lambda( \partial^\beta \sigma) 
\frac{\up^\beta}{\beta!} +\sum_{\vert \beta\vert=N+1}\frac{\up^\beta }{\beta!}  \, \int_0^1 (1-t)^N\, \rho \left(  \partial^\beta\sigma (\cdot+
  t\up)\right)\,   dt=0$$
with respect to the coordinates of $\up$ at $\up=0$ yields the first part of
the assertion.
\item 
Let $\lambda_1$ and $\lambda_2$ be two $\Z^d$-translation invariant
linear forms on $ CS^{\notin \Z}(\R^d)$ satisfying the assumptions of the
theorem with the same restriction $\rho$.  The   Taylor formula
(\ref{eq:Taylorexpansion}) applied to $\sigma \in CS_{\rm c.c}^{\notin \Z}(\R^d)$  yields by linearity  of $\lambda_i$ and for $N$ chosen large enough: $$\lambda_i\left( t_{\up}^* \sigma\right)=\sum_{\vert \beta\vert\leq N} \lambda_i \left(
 \partial^\beta \sigma \right)
\frac{\up^\beta}{\beta!} +\sum_{\vert \beta\vert=N+1}\frac{\up^\beta }{\beta!} \, \int_0^1 (1-t)^N\, \rho \left( \partial^\beta\sigma (\cdot+
  t\up)\right)\,   dt,$$ 
so that $\  t_{\up}^*\lambda_1(\sigma)-\sum_{\vert \beta\vert\leq 
  N}\lambda_1(\partial^\beta\sigma)\, \frac{\up^\beta}{\beta!} = t_{\up}^*\lambda_2(\sigma)-\sum_{\vert \beta\vert\leq 
  N}\lambda_2(\partial^\beta\sigma)\, \frac{\up^\beta}{\beta!} \quad \forall  \up \in \Z^d$ since 
$\lambda_1$ and $\lambda_2$ both coincide with $\rho$   on $CS_{\rm
  c.c}^{<-d}(\R^d)\cap CS_{\rm c.c}^{\notin \Z}(\R^d)$. Since
$t_{\up}^*\lambda_i= \lambda_i$ for any $\up \in \Z^d$,
this implies that the
polynomial expressions $\sum_{0<\vert \beta\vert\leq 
  N}\lambda_1(\partial^\beta\sigma)\, \frac{\up^\beta}{\beta!}$ and $\sum_{0<\vert \beta\vert\leq 
  N}\lambda_2(\partial^\beta\sigma)\, \frac{\up^\beta}{\beta!}$ in the coordinates
of $\up$ coincide for all $\up \in \Z^d$  and hence that  their coefficients
coincide $\lambda_1(\partial^\beta\sigma)=\lambda_2(\partial^\beta\sigma)$ when $0<
\vert \beta\vert <N$. Since this holds for any large enough $N$, we conclude
that $\lambda_1(\partial^\beta\sigma)=\lambda_2(\partial^\beta\sigma)$
when $\beta\neq 0$. It follows that its value on derivatives  $\lambda(\partial^\beta \sigma),\quad  \beta \neq 0, \sigma\in CS_{\rm c.c}^{\notin \Z}(\R^d)$  is entirely
  determined by the restriction $\rho$ to $CS_{\rm c.c}^{<-d}(\R^d)$.
\end{enumerate}\endsquare
\subsection{The canonical regularised integral and discrete sum}In order to determine all translation invariant linear forms on $ CS_{\rm
  c.c}^{\notin \Z}(\R^d)$ that extend the integral and discrete sum on
$L^1(\R^d)$ we   need the following lemma which collects results from \cite{FGLS} (see
Lemma 1.3), see also \cite{P2}.\begin{lem} \label{lem:FGLS}
Any  symbol $\sigma\in CS^{\notin\Z}_{\rm c.c}(\R^d)$ is up to some
smoothing symbol,
  a finite sum  of partial
  derivatives, i.e. there exist symbols  $\tau_i\in CS^{\notin \Z}_{\rm c.c}(\R^d), i=1,
  \cdots, n$  such that
 \begin{equation}\label{eqsigmapartial}\sigma \sim \sum_{i=1}^n
   \partial_i\tau_i. \end{equation}
\end{lem}
 The following theorem shows that translation
invariant linear forms on non integer order symbols are entirely determined by
their restriction to $L^1$-symbols. 
\begin{thm}\label{thm:uniqueness}  Let $\lambda: CS_{\rm c.c}^{\notin \Z}(\R^d)\to \C$ be a
  linear form   which coincides on
  $CS_{\rm c.c}^{\notin \Z}(\R^d)\cap CS_{\rm c.c}^{<-n}(\R^d) $  with a
  linear form $\rho: L^1(\R^d)\to \C $ satisfying (\ref{eq:rho}). \\
 If $\lambda$ is either $\Z^d$ or $\R^d$-translation invariant, then it
  is entirely determined by its restriction $\rho$.
\end{thm}
\begin{rk} $\R^d$-translation invariance and its relation to Stokes' property
  had already been investigated in \cite{P2}. 
\end{rk}   
{\bf Proof:}  Let $\lambda$ be an $\R^d$ or $\Z^d$-translation invariant
linear form on $ CS^{\notin \Z}(\R^d)$ satisfying the asumptions of the
theorem. Since a symbol  $\sigma\in CS_{\rm
  c.c}^{\notin \Z}(\R^d)$ has vanishing residue by Lemma \ref{lem:FGLS},   up to a
smoothing symbol, it is  a finite linear
combination of partial derivatives of symbols in $CS_{\rm
  c.c}^{\notin \Z}(\R^d)$. With the notations of the lemma we write 
$\sigma= \sum_{i=1}^n \partial_i \tau_i+r$ for some smoothing symbol
  $r $ on which by assumption, the linear form $\lambda$ coincides with its
  restriction $\rho$
   to
$CS^{<-n}(\R^d)$. Thus, by linearity
$$\lambda(\sigma )=\sum_{i=1}^n
   \lambda(\partial_i\tau_i)+\rho(r). $$
But by Proposition \ref{prop:extrho} the value of $\lambda$ on derivatives is entirely
determinaed by its restriction $\rho$. It follows that $\lambda$ on
$CS^{\notin \Z}(\R^d)$ is entirely determined by its restriction $\rho$ as announced in the theorem.\endsquare
\\ \\
On the  grounds of Theorem \ref{thm:uniqueness}, we define the canonical
regularised integral over $\R^d$ and the canonical regularised sum over $\Z^d$ on non integer order symbols.
\begin{defn} Let the canonical regularised integral $\cutoffint_{\R^d}$ (resp.
  the canonical regularised
sum $\cutoffsum_{\Z^d}$) be the unique $\R^d$ (resp. $\Z^d$) -translation
invariant extension of the Riemann integral $\int_{\R^d}$ (resp. sum
$\sum_{\Z^d}$) to $CS_{\rm c.c}^{\notin \Z}(\R^d)$. 
\end{defn} 
\vfill \eject \noindent
\section{Concrete realisations  of the canonical regularised integral and discrete sum}
\subsection{Integrals of symbols over infinite convex polytopes }
 We recall from \cite{P2} that the  Hadamard finite part integral $$\sigma\mapsto {\rm fp}_{R\to
    \infty}\int_{B(0, R)}\sigma(\ux)\, d\ux $$ defines an $\R^d$-translation
  invariant linear form on $CS_{\rm c.c}^{\notin \Z}(\R^d)$ which extends the Riemann
  integral on $L^1$ symbols. \\ By Theorem \ref{thm:uniqueness} it therefore coincides with the canonical
  integral $$\cutoffint_{\R^n}\sigma= {\rm fp}_{R\to
    \infty}\int_{B(0, R)}\sigma(\ux)\, d\ux\quad\forall \sigma\in CS_{\rm c.c}^{\notin
    \Z}(\R^d).$$
\begin{rk} For general symbol $\sigma\in CS_{\rm c.c}(\R^d)$, the expression  ${\rm fp}_{R\to
    \infty}\int_{B(0, R)}\sigma(\ux)\, d\ux$ depends on the choice of the
  parameter $R$ but for convenience and following the notations of \cite{L} we  extend this notation to all symbols setting  $\cutoffint_{\R^n}\sigma:= {\rm fp}_{R\to
    \infty}\int_{B(0, R)}\sigma(\ux)\, d\ux$ for any $\sigma\in CS_{\rm c.c}(\R^d)$,
  which we refer to as the cut-off regularised integral of $\sigma$ as recalled
  in formula (\ref{eq:constanttermclassical}) of the Appendix.
\end{rk} 
In the following, we show that the canonical regularised integral on non
integer order symbols  can also be realised as a finite part as $R$
tends to infinity of an integral over an $R$-expanded polytope. \\
Let us first introduce some notations, following those of \cite{GSW}, see also
\cite{KSW2}. \\
A compact convex polytope  $\Delta$ in $\R^d$ is a compact set which can be
obtained as the intersection of finitely many half-spaces $\Delta=H_1\cap\cdots
\cap H_m$ where 
$$H_i=\{x\in \R^d, \langle u_i, x\rangle +a_i\geq 0\}$$  for some vectors
$u_i\in \R^d,
i=1, \cdots, m$ and some positive real numbers $a_i, i=1, \cdots, m$. \\
 Assuming that
  $m$ is the smallest possible number of such
half spaces, the {\it facets } or codimension one faces of $\Delta$ are the $\sigma_i=
\Delta\cap \partial H_i, i=1, \cdots, m$. The vector $u_i$ can be thought of
 as an  inward normal vector to the $i$-th
facet $\sigma_i$.\\
Alternatively, a compact convex
polytope is the convex hull of a finite set of points in $\R^d$ and when this
set is minimal, its elements   are called
{\it vertices } of $\Delta$. \\
A compact convex polytope   is {\it simple} if each vertex is the intersection of
exactly $d$ facets. It is called {\it integral} or {\it lattice polytope} if its
vertices lie in the lattice $\Z^d$. It is {\it regular} if moreover, the edges
emanating for each vertex lie along lines generated by a $\Z$-basis of the
lattice $\Z^d$ or equivalently if   each local cone at any of its vertices can be
 transformed by an integral unimodular affine transformation to a
 neighborhood of $0$. \\ 
 Let  $\Delta\subset \R^d$ be a simple, integral and regular compact convex
 polytope and such that the origin $0$ lies in the inerior of $\Delta$. Since
 $\Delta$ is integral, the vectors $u_i, i=1, \cdots, m$ can be chosen to belong
 to the dual lattice $\Z^{d*}$  and they can be normalized to be {\it primitive
 lattice elements} i.e. not be a multiple of a lattice element by an integer
 larger than one. As pointed out in \cite{KSW2}, the fact that a normal vector
 $u$ to a facet $\sigma$ can be chosen to be integral is a consequence of
 Cramer's rule; indeed, choosing integral edge vectors  $\beta_1, \cdots,
 \beta_d$ that emanate from a vertex on the facet $\sigma$ and such that
 $\beta_1, \cdots, \beta_{d-1}$ span the tangent plane to the facet and
 $\beta_n$ is transverse to $\sigma$, solving the linear equations $\langle u,
 \beta_1\rangle=\cdots= \langle u, \beta_{d-1}\rangle=0$ and $\langle u,
 \beta_d\rangle=1$, we get an inward normal vector $u$ with rational entries,
 from which it is easy to build an integral inward normal vector. \\
 The expanded polytope 
$R\cdot \Delta$  for any positive real number $R$ is defined by:
$$\ux\cdot u_i+R a_i\geq 0, \quad i=1, \cdots ,m.$$ 
Note that $\bigcup_{N\in \N} N\cdot \Delta= \R^d.$
 \begin{ex}Let $\vert \ux\vert_{\rm sup}:= {\rm sup}_{i=1}^d \vert x_i\vert$ denote the
supremum norm of $\ux=(x_1, \cdots, x_d)\in \R^d$. 
 The balls $B_{\rm sup}(0, R)=\{\ux\in \R^d, \vert \ux\vert_{\rm sup } \leq
 R\}$  for the supremum norm  are  defined by the inequalities: 
$$-R\leq \langle\ux, e_i\rangle\leq R, \quad i=1, \cdots ,d$$
where $(e_1, \cdots, e_d)$ is the canonical orthonormal basis in $\R^d$. 
Setting  $m=2d$, $u_i=e_i, i=1, \cdots, d$,  $a_i=1, i=1, \cdots, d$, 
$u_{i+d}=-e_i,i=1, \cdots, d$, $a_{i+d}=1, i=1, \cdots, d$, we recover the
above set of inequalities which define an expanded  simple  integral regular
(compact convex) polytope.
\end{ex}
For a continuous  function
$f$ on $\R^d$                and any $R>0$ we  set:
$$\tilde P_{R\cdot \Delta}(f):=\int_{R\cdot \Delta}f(\ux)\, d\ux.$$

\begin{prop}\label{prop:intpolytope} Let $\sigma \in CS^a_{\rm c.c}(\R^d)$ for
  some complex number $a$.\begin{enumerate}
\item The map 
$R\mapsto \tilde  P_{R\cdot \Delta}(\sigma)$  defines a log-polyhomogeneous symbol of
order $a+d$  and log-type $1$ on $\R^+$ which differs from the
log-polyhomogenous symbol  $R\mapsto
\int_{B(0, R)}\sigma(\ux)\, d\ux$ by a classical symbol of order $a+d$.
\item  The
constant term ${\rm fp}_{R\to \infty}\tilde  P_{R\cdot \Delta}(\sigma)  $ in its
 asymptotic expansion relates to  the finite part ${\rm fp}_{R\to \infty}\cutoffint_{B(0, R)}\sigma(\ux)\,
d\ux$ as follows
$${\rm fp}_{R\to \infty}\tilde  P_{R\cdot \Delta}(\sigma)-{\rm fp}_{R\to \infty}\int_{B(0, R)}\sigma(\ux)\,
d\ux 
= \int_{\Delta\ominus B(0, 1)} \sigma_{-d}(\ux)\,
d\ux\quad\forall \sigma\in CS_{\rm
  c.c}(\R^d),$$
with $A\ominus B:= (A-(A\cap B))\cup (B-(A\cap B))$ the symmetric
difference.\footnote{We do not use the usual notation $A\Delta B$ to avoid a
  clash of notations with the polytope commonly denoted in the literature by
  $\Delta$.}
\end{enumerate}
\end{prop}
\begin{rk} In particular for any $\sigma \in CS_{\rm c.c}(\R^d)$
  the map $R\mapsto \int_{\vert \ux\vert_{\rm sup}\leq R} \sigma(\ux)\,
d\ux$ defines  a log-polyhomogeneous symbol of
order $a+d$  and log-type $1$ and  if $\sigma$ has non integer order 
  we have
\begin{equation}\label{eq:supnormint}{\rm fp}_{R\to \infty}\int_{\vert \ux\vert_{\rm sup}\leq R} \sigma(\ux)\,
d\ux= {\rm fp}_{R\to \infty}\int_{B(0, R)} \sigma(\ux)\, d\ux\quad \forall \sigma \in CS_{\rm
  c.c}^{\notin \,\Z}(\R^d).
\end{equation}
\end{rk}
{\bf Proof:} By  formula (\ref{eq:Rasymptotics}) of the Appendix,
   the map $R\mapsto \int_{B(0, R)} \sigma(\ux)\, d\ux$  defines a log-polyhomogeneous symbol of order
$a+d$ and log-type $1$.
 Let us  estimate the difference  for large $R$:
\begin{eqnarray}\label{eq:intdeltaB}
&{}&\int_{R\, \Delta} \sigma(\ux)\, d\ux -\int_{B(0, R)} \sigma(\ux)\, d\ux \nonumber\\
&=&   \sum_{j=0}^{N} \int_{R\, \Delta\ominus B(0, R)}\chi(\ux) \sigma_{a-j}(\ux)
\, d\ux+\int_{R\, \Delta\ominus B(0, R)}\sigma_{(N)}(\ux)\, d\ux\nonumber\\
&=& \sum_{j=0}^{N} \int_{R\, \Delta\ominus B(0, R)} \sigma_{a-j}(\ux)\, d\ux
+\int_{R\, \Delta\ominus B(0, R)}\sigma_{(N)}(\ux)\, d\ux,
\end{eqnarray} where we
use   the fact that
$\chi$ is one outside the unit euclidean ball.\\ Since
$\sigma_{(N)}$ is a symbol of order $a-N$ we have:
$$\vert \sigma_{(N)}(\ux)\vert\leq C(1+\vert \ux\vert^2)^{\frac{{\rm Re}(a)-N}{2}} $$
for some constant $C$. Hence, for large enough $N$
\begin{eqnarray*}
 \vert \int_{R\, \Delta- (R\Delta\cap B(0, R))}\sigma_{(N)}(\ux)\,
d\ux\vert
&\leq &C\, (1+R^2)^{\frac{{\rm Re}(a)-N}{2}}{\rm Vol}\left(R\, \Delta-
  (R\Delta\cap B(0, R))\right)\\
&\leq &C\, R^d(1+R^2)^{\frac{{\rm Re}(a)-N}{2}}{\rm Vol}\left(
  \Delta\right)\\
&\leq &C\, (1+ R^2)^{\frac{{\rm Re}(a)+d-N}{2}}{\rm Vol}\left( \Delta\right).
\end{eqnarray*}
Using once more the fact that 
$\sigma_{(N)}$ is a symbol of order $a-N$ combined with the equivlence of the
supremum and the euclidean norms we also have:
$$\vert \sigma_{(N)}(\ux)\vert\leq C^\prime\, (1+\vert \ux\vert_{\rm sup})^{{\rm Re}(a)-N} $$ for some constant $C^\prime$
and
\begin{eqnarray*}
 \vert \int_{B(0, R)- (R\Delta\cap B(0, R))}\sigma_{(N)}(\ux)\,
d\ux\vert
&\leq &C^\prime\, (1+d\, R)^{{\rm Re}(a)-N}{\rm Vol}\left(B(0, R)-
  (R\Delta\cap B(0, R))\right)\\
&\leq &C^\prime\, R^d(1+R)^{{\rm Re}(a)-N}{\rm Vol}\left(B(0, 1)\right)\\
&\leq &C^\prime\, (1+ R)^{{\rm Re}(a)+d-N}{\rm Vol}\left( B(0, 1)\right).
\end{eqnarray*}
Consequently,  we can  choose 
$N$ sufficiently  large  so that
$$
\int_{R\, \Delta\ominus B(0, R)}\sigma_{(N)}(\ux)\, d\ux
= O((1+R)^{{\rm Re}(a)+d-N}).$$
This settles the case of  integrals involving the  remainder term $\sigma_{(N)}$. As for  integrals of homogeneous symbols $\int_{R\, \Delta\ominus B(0, R)}
\sigma_{a-j}(\ux)\, d\ux$, we have$$
\int_{R\, \Delta \ominus B(0, R) } \sigma_{a-j}(\ux)\, d\ux 
=  
\int_{\Delta\ominus B(0, 1) } \sigma_{a-j}(R \, \uy)\, R^d\,  d\uy 
= R^{a-j+d}
\int_{\Delta\ominus  B(0, 1) } \sigma_{a-j}( \uy)\,  d\uy, $$
which shows they are homogeneous of degree $a-j+d$. \\
Combining these results shows that $R\mapsto \int_{R\, \Delta \ominus B(0, R)
}\sigma(\ux)\, d\ux$ defines a classical symbol of order $a+d$
with constant term given by 
$${\rm fp}_{R\to \infty}\int_{R\, \Delta \ominus B(0, R)
}\sigma(\ux)\, d\ux=
\int_{\Delta\ominus  B(0, 1) } \sigma_{-d}( \uy)\,  d\uy. $$
Since $R\mapsto \int_{ B(0, R)
}\sigma(\ux)\, d\ux$  defines a log-polyhomogeneous symbol of order
$a+d$ and log-type $1$, so does the map $R\mapsto \int_{ R\, \Delta
}\sigma(\ux)\, d\ux$  by (\ref{eq:intdeltaB}). Its finite part differs from
${\rm fp}_{R\to \infty} 
 \int_{ B(0, R)
}\sigma(\ux)\, d\ux$ by ${\rm fp}_{R\to \infty}\int_{R\, \Delta \ominus B(0, R)
}\sigma(\ux)\, d\ux=
\int_{\Delta\ominus  B(0, 1) } \sigma_{-d}( \uy)\,  d\uy. $ \endsquare
\subsection{Discrete sums of classical symbols on infinite convex polytopes} 
With  notations similar to those   of \cite{GSW} and \cite{KSW2},  for any  $R>0$ we define 
the {\it dilated} polytope $\Delta_{R,h}\subset \R^d$ obtained by shifting the $i$-th facet
of $\Delta_R$ outwards by a distance $h_i$. It is described by the inequalities:
$$ \langle \ux,u_i\rangle+R\, a_{i}+h_{i}\, \geq 0,\quad i=1, \cdots, m$$ with $h=(h_1,
\cdots, h_m)\in \R^m$. \\
When $\Delta$ is simple, so is $\Delta(h)$ simple for $h$ small enough. 
For $h=0$, $\Delta_{R, 0}$  yields back the expanded polytope  $R\, \Delta$ and for $R=0$ we set
  $\Delta_{h}:= \Delta_{0,h}$.
For a function
$f$ on $\R^d$                and any $R>0$,  any $ h\in
\R^{2d}$ we further set:
$$ P_{R\, \Delta}(f)( h):=\sum_{\Delta_{R, h}\cap \Z^d}f(\un),\quad \tilde P_{R\, \Delta}(f)( h):=\int_{\Delta_{R, h}}f(\ux)\, d\ux.$$
For $h=0$, the integral, resp. the sum are taken over the expanded polytope $R\, \Delta$,
resp. integer points in $R\, \Delta$ and we set 
$P_{R\, \Delta}(f):= P_{R\, \Delta}(f)(0)$ similarly to 
 $\tilde P_{R\, \Delta}(f)=\tilde P_{R\, \Delta}(f)( 0)$.\\ \\
 Let us first establish a technical lemma which will be useful for what follows.
\begin{lem}\label{lem:technicallem}Given a symbol $\sigma\in CS^a_{\rm
    c.c}(\R^d)$ with complex order $a$ and any positive
  integer $\gamma_i$,  the map
 $R\mapsto \partial^{\gamma_i}_{h_i}\left(\tilde P_{R\, \Delta}(\sigma)( h)\right)_{\vert_{h=0}}$ defines a log-polyhomogeneous symbol on $\R^+$ of order
$a+d-\gamma_i$. Its finite part given by the constant term in the symbolic
asymptotic expansion in $R$ as $R\to \infty$ reads:
$${\rm fp}_{R\to \infty} \left( \partial^{\gamma_i}_{h_i}\tilde P_{R\, \Delta}(\sigma)( h)\right)_{\vert_{h=0}}=
\left(\partial^{\gamma_i}_{h_i}\tilde P_{
    \Delta}(\sigma_{-d+\gamma_i-1})(h)\right)_{\vert_{h=0}}$$ and therefore  vanishes whenever the order of $\sigma$ is non integer.
\end{lem}
\begin{rk} The case $\gamma_i=0$ was dealt with in
  Proposition  \ref{prop:intpolytope} where we showed that $R\mapsto  \tilde P_{R\, \Delta}(\sigma)$
defines a log-polyhomogeneous symbol of order $a+d$ and  log-type $1$. 
\end{rk}
{\bf Proof:} \begin{enumerate}
\item Let us  first carry out the proof in the case of a hypercube  defined by
the following set of inequalities 
$$-\beta_i\leq \langle \ux, e_i \rangle\leq \alpha_i, \quad i=1, \cdots ,d$$ 
for given  positive numbers  $\alpha_1, \cdots,
\alpha_d, \beta_1, \cdots, \beta_d$. Here as before, $\{e_1, \cdots, e_d\}$ stands for the canonical
orhtonormal basis of $\R^d$. Setting $m=2d$, we rewrite the inequalities
defining $\Delta$ as follows:
$$ \langle \ux, u_i\rangle+a_{i}\, \geq 0; \quad \langle \ux, u_{i+d}\rangle+ a_{i+d}\, \geq 0,$$
where we have set $u_i=-e_i,a_i=\alpha_i,  u_{i+d}= e_i, a_{i+d}=\beta_i$.
\\For any positive real number $R $,  
 \begin{eqnarray*}
&{}&\left(\partial^{\gamma_i}_{h_i}\tilde P_{R\, \Delta}(\sigma)(h)\right)_{\vert_{h=0}}\\
&=& \left(\int_{-h_{1}-R\, \beta_1\leq x_1\leq
    h_{1+d}+R\, \alpha_1}dx_1 \cdots \partial^{\gamma_i}_{h_i}
  \int_{-h_{i}-R\, \beta_i\leq x_i\leq
    h_{i+d}+R\, \alpha_i}dx_i \cdots \right.\\
&{}&\left. \cdots    \int_{-h_{d}-R\, \eta_d\leq x_d\leq h_{2d} +R \, \alpha_d} dx_d\, \sigma(\ux)\right)_{\vert_{h=0}}\\
&=&(-1)^{\gamma_i}\prod_{j=1, j\neq i}^d
\int_{-R\, \beta_j\leq x_j\leq R\, \alpha_j}dx_j  \, \partial_i^{\gamma_i-1}\sigma\left(x_1,
  \cdots,x_{i-1},-R\, \beta_i, x_{i+1}, \cdots, x_d
  \right)\\
\end{eqnarray*}
is the integral over the face  $x_i=-R\beta_i$ of the expanded polytope $R\,
\Delta$. \\
 If $\sigma\in CS_{\rm c.c}^a(\R^d)$ then $\tau_i:=\partial^{\gamma_i-1 }\sigma$
 defines a classical symbol of order $a_i=a-\gamma_i+1$   on  $
 \R^{d}$. With  the notations of (\ref{eq:asymptsymb}) in the Appendix, let
 us write
\begin{equation}\label{eq:useful3} \tau_i(\ux) =\sum_{k=0}^N\chi(\ux)\,\tau_{i, a_i-k}(\ux) +\tau_{i,
  N}(\ux).
\end{equation} Provided $R$
is chosen large enough, on the face $x_i=-R\beta_i$ of the expanded polytope we have
\begin{eqnarray*}
&{}&\tau_i\left(x_1,
  \cdots,x_{i-1},-R\, \beta_i, x_{i+1}, \cdots, x_d
  \right) \\
&=&\sum_{k=0}^N\tau_{i, a_i-k}\left(x_1,
  \cdots,x_{i-1},-R\, \beta_i, x_{i+1}, \cdots, x_d
  \right) +\tau_{i,
  N}\left(x_1,
  \cdots,x_{i-1},-R\, \beta_i, x_{i+1}, \cdots, x_d
  \right).
\end{eqnarray*}Hence,
\begin{eqnarray}\label{eq:useful1}
&{}&\prod_{j=1, j\neq i}^d
\int_{-R\, \beta_i\leq x_j\leq R\, \alpha_i}dx_j  \, \partial_i^{\gamma_i-1}\sigma\left(x_1,
  \cdots,x_{i-1},-R\, \beta_i, x_{i+1}, \cdots, x_d
  \right)\nonumber\\
&=&\sum_{k=0}^N 
R^{a_i+d-k}\, \prod_{j=1, j\neq i}^d\int_{-\, \beta_j\leq x_j\leq \, \alpha_j}dx_j  \, \tau_{i,a_i-k}\left(x_1,
  \cdots,x_{i-1},-\beta_i, x_{i+1}, \cdots, x_d
  \right)\nonumber\\
&+&\prod_{j=1, j\neq i}^d
\int_{-R\, \beta_j\leq x_j\leq R\, \alpha_j}dx_j \tau_{i, N}\left(x_1,
  \cdots,x_{i-1},-R\, \beta_i, x_{i+1}, \cdots, x_d
  \right). 
\end{eqnarray}
Since $\tau_{i, N}$ is a symbol of order $a_i-N$, we have
$$\vert\tau_{i, N}(\ux)\vert\leq C\, \left(1+\vert
  \ux\vert^2\right)^{\frac{{\rm Re}(a_i)-N}{2}}$$
for some constant $C$.  On the face of the polytope defined by equation
$x_i=-R\, \beta_i$,  choosing $N$ large enough for ${\rm Re}(a_i)-N$ to be negative
we write
$$\vert\tau_{i, N}(\ux)\vert\leq C\, \left(1+\beta_i^2 \,
  R^2\right)^{\frac{{\rm Re}(a_i)-N}{2}}.$$
It follows that  there is some large enough $N$
such that  
 \begin{eqnarray*}
&{}&\vert \prod_{j=1, j\neq i}^d
\int_{-R\, \beta_j\leq x_j\leq R\, a_j}dx_j \tau_{i, N}\left(x_1,
  \cdots,x_{i-1},-R\, \beta_i, x_{i+1}, \cdots, x_d
  \right)\vert\\
 &\leq& C\, R^d\,(1+\beta_i^2 \,
  R^2)^{\frac{{\rm Re}(a_i)-N}{2}}\prod_{j=1, j\neq i}^d
\int_{-\, \beta_j\leq x_j\leq \, a_j}dx_j\\
&\leq& C^\prime (1+ 
  R^2)^{\frac{{\rm Re}(a_i)+d-N}{2}} 
\end{eqnarray*}
 for some constant $C^\prime$.\\
Combining these results shows that the map  $R\mapsto \left(\partial^{\gamma_i}_{h_i}\tilde P_{R\,
    \Delta}(\sigma)(h)\right)_{\vert_{h=0}}$ defines a classical symbol of
order $a_i+d=a+d-\gamma_i+1$. \\ Using (\ref{eq:useful1}) we see that  its 
finite part  as $R$ tends to infinity reads
\begin{eqnarray*}
&{}&{\rm fp}_{R\to \infty}\left(\partial^{\gamma_i}_{h_i}\tilde P_{R\,
    \Delta}(\sigma)(h)\right)_{\vert_{h=0}}\\
&=& 
 \prod_{j=1, j\neq i}^d\int_{-\, \beta_j\leq x_j\leq \, \alpha_j}dx_j  \,\tau_{i,-d}\left(x_1,
  \cdots,x_{i-1},-\, \beta_i, x_{i+1}, \cdots, x_d
  \right)\\
&=&
 \prod_{j=1, j\neq i}^d\int_{-\, \beta_j \leq x_j\leq \, \alpha_j}dx_j  \,\left( \partial^{\gamma_i-1}\sigma\right)_{-d}\left(x_1,
  \cdots,x_{i-1},-\, \beta_i, x_{i+1}, \cdots, x_d\right)\\
&=& \left(\partial^{\gamma_i}_{h_i}\tilde P_{  \Delta}(\sigma_{-d+\gamma_i-1})(h)\right)_{\vert_{h=0}}\\
\end{eqnarray*}
and therefore vanishes if $\sigma$  has non integer order.
\item The proof generalises to  any general simple regular integral convex
  polytope. Indeed, using (\ref{eq:useful3}), we split as before the
  derivative $\left(\partial_{h_i}^{\gamma_i}\tilde P_{R\Delta}(\sigma)(h)\right)_{\vert_{h=0}}$ into a
    finite sum of homogeneous terms and a remainder term. The   faces of the
    polytope moving out to infinity as
    $R\to \infty$, the  behaviour as
    $R\to \infty$ is controlled using 
\begin{itemize}\item on the one hand, the homogeneity of the derivatives of the integrals of the
    homogenous terms which turn out to be integrals of derivatives of the
    homogeneous components of the 
    symbol over   faces of the
    polytope and 
\item on the other hand, the symbolic
    behaviour of the integrand in the remainder term which is an integral over
    faces of the polytope of a symbol the  order of which has arbitrarily
    negative real part.
\end{itemize} 
\end{enumerate}\endsquare\\ \\
We now want to compare $P_{N\,\Delta}(\sigma)$ and $\tilde P_{N\,\Delta}(\sigma)$
as $N$ goes to infinity.\\
Following the notations of \cite{GSW}, let 
$${\rm Todd}\left(\partial\right):= \sum_\alpha b_\alpha \partial^\alpha$$
where $\sum_\alpha b_\alpha \partial^ \alpha $ is the Taylor series expansion at the origin of the Todd function: 
$${\rm Todd}(x)= \prod_{i=1}^m \frac{x_i}{1-e^{-x_i}}.$$
The Khovanskii-Pukhlikov formula \cite{KP} compares the discrete sum  $P_{N\, \Delta}(f)$ with the
integral $\tilde P_{N\, \Delta}(f)$ over the expanded convex polytope  $N\, \Delta$ for polynomial functions $f$:
$$  P_{N\, \Delta}(f) -\tilde P_{N\, \Delta}(f)=\left(\left( {\rm Todd}(\partial_h)-Id\right)\tilde P_{N\, \Delta}(f)(h)\right)\vert_{h=0}.$$
The following proposition is a reformulation of results of  \cite{GSW} (formula (15), see also
\cite{AW},  \cite{KSW1} and 
\cite{KSW2} for previous results along these lines) where the authors 
generalise the
Khovanskii-Pukhlikov formula to classical symbols, 
in which  case the formula is
not exact anymore but only holds asymptotically. 
\begin{prop}\label{prop:GSW} \cite{GSW}  Given a symbol $\sigma\in CS^a_{\rm
    c.c}(\R^d)$ with complex order $a$, the discrete map $N\mapsto P_{N\, \Delta}(\sigma)$ can be interpolated
  by a log-polyhomogeneous symbol on $\R^+$  of order $a+d$ and log-type $1$:
\begin{equation}\label{eq:interpolation} R\mapsto \tilde P_{R\, \Delta}(\sigma)+ \left( {\rm
    Todd}(\partial_h)-Id\right) \tilde  P_{R\, \Delta}(\sigma) (
h)\vert_{h=0}+C(\sigma)
\end{equation} for some constant $C(\sigma)$ independent of the
choice of expanded polytope, i.e.
$$P_{N\, \Delta}(\sigma)- \tilde P_{N\, \Delta}(\sigma)\sim_{N\to \infty} \left( {\rm
    Todd}(\partial_h)-Id\right) \tilde  P_{N\, \Delta}(\sigma) (
h)\vert_{h=0}+C(\sigma).$$
More precisely, there are polynomials $M^{[j]},j\in \N$ on $\R^d$ such that 
\begin{equation}\label{eq:GSW}
 P_{N\, \Delta}(\sigma) -\tilde P_{N\, \Delta}(\sigma) = \left(\left(
   M^{[j]}(\partial_h)-Id\right) \tilde P_{N\, \Delta}(\sigma) (
 h)\right)_{\vert_{h=0}}+R^j(\sigma)( N)
\end{equation}
where $$R^j(\sigma)( N):= \sum_{p}(-1)^p \int_{C_{p, N} }\sum_{\vert
  \alpha\vert=j}^{\vert \alpha\vert=d\,j} \phi^p_{\alpha, j}
(\ux)\,\partial^\alpha \sigma( \ux)\, d\ux$$ tends to $C(\sigma)$ as $N\to \infty$.\\
The $C_{p, N}$ are convex polytopes  growing with $N$ and $\phi^p_{\alpha,j}$
bounded piecewise smooth periodic functions as described in \cite{KSW2} and \cite{AW}.
\end{prop}
\begin{rk} Let us comment on the interpolation (\ref{eq:interpolation}) which is  a reinterpretation
of the statements of \cite{GSW}. \\
Since we know by Lemma \ref{lem:technicallem} that for any positive $\gamma_i$,   the maps $R\mapsto \left(
  \partial_{h_i}^{\gamma_i} \tilde P_{R\, \Delta}(\sigma) (
 h)\right)_{\vert_{h=0}}$ define classical  symbols on $\R^+$ of
order $a+d-\gamma_i$  if $\sigma$
 is of order $a$, it follows that  the maps $R\mapsto \left(\left(
   M^{[j]}(\partial_h)-Id\right) \tilde P_{R\, \Delta}(\sigma) (
 h)\right)_{\vert_{h=0}}$  give rise to classical  symbols on $\R^+$   of order
 $a+d-k_j$ for some integer $k_j$ and log-type $1$. This combined with
 (\ref{eq:GSW})  tells us that  (\ref{eq:interpolation}) indeed provides an
 interpolation of $N\mapsto P_{N\, \Delta}(\sigma)$ by a log-polyhomogeneous
 symbol of order $a+d-\gamma_i$ and log-type $1$.  \end{rk}
Since  the expression  $P_{N\, \Delta} (\sigma)=\tilde P_{N\, \Delta}(\sigma) + \left(\left(
   M^{[j]}(\partial_h)-Id\right) \tilde P_{N\, \Delta}(\sigma) (
 h)\right)_{\vert_{h=0}}+R^j(\sigma)( N)$ has an asympotic expansion as $N\to
\infty$  of the same typs as  $\tilde P_{N\, \Delta}(\sigma)$, we can extract
the constant term in the expansion and consier the Hadamard finite part $ {\rm fp}_{N\to \infty} \sum_{N\, \Delta\cap
   \Z^d}\sigma(\un).$
\begin{rk} Clearly, when $\sigma \in CS^{<-d}(\R^d)$ then  $ {\rm fp}_{N\to \infty} \sum_{N\, \Delta\cap
   \Z^d}\sigma(\un)= \lim_{N\to \infty} \sum_{N\, \Delta\cap
   \Z^d}\sigma(\un)= \sum_{ \Z^d}\sigma(\un)$ is an ordinary limit.
\end{rk}
\begin{ex} Given a polynomial $Q(\ux)=  \sum_{\vert \alpha\vert \leq K} a_\alpha
  \ux^\alpha$ on $\R^d$, $a_\alpha\in \C$  then:
$${\rm fp}_{N\to \infty} \cutoffsum_{\vert \un\vert_{\rm sup}\leq N} Q(\ux)= 0. $$
\\ Indeed, we have 
$$\sum_{\vert \un\vert_{\rm sup}\leq N} Q(\un)= \sum_{\vert \alpha\vert \leq
  K} a_\alpha\prod_{i=1}^d 
\sum_{-N}^N n_i^{\alpha_i}.$$
Since the one dimensional sum $N\mapsto \sum_{-N}^N n_i^{\alpha_i}$ is known
to be  polynomial in $N$ (this easily follows from the Euler-MacLaurin
formula, see e.g. \cite{MP}), its finite part as $N\to \infty$ coincides with its
value at $N=0$: 
\begin{eqnarray*}
{\rm fp}_{N\to \infty}\sum_{\vert \un\vert_{\rm sup}\leq N} Q(\un)
&=&  {\rm fp}_{N=0}\sum_{\vert \un\vert_{\rm sup}\leq N} Q(\un)\\
&=& \sum_{\vert \alpha\vert \leq K} a_\alpha\prod_{i=1}^d
{\rm fp}_{N=0} \sum_{-N}^N n_i^{\alpha_i}.
\end{eqnarray*}
One easily checks that this last expression ${\rm fp}_{N=0} \sum_{-N}^N
n_i^{\alpha_i}$  vanishes  as a result of the fact that    ${\rm fp}_{N=0}
\sum_{n=-N}^N P(n)=0$ for any polynomial $P$ \cite{MP}.
 \end{ex}
\begin{thm} \label{thm:translinvregsum}
The  Hadamard finite part  $\sigma\mapsto {\rm fp}_{N\to \infty}\sum_{\vert \un\vert_{\rm sup}\leq N}
\sigma(\un)$ coincides with the canonical
regularised sum on non integer order symbols and we have:
\begin{equation} \label{eq:cutoffsumnonintorder}{\rm fp}_{N\to \infty} \sum_{N\,
  \Delta\cap \Z^d} \sigma= \cutoffsum_{\Z^d} \sigma=\cutoffint_{\R^d}\sigma+C(\sigma)\quad \forall
\sigma\in CS_{\rm c.c}^{\notin\Z}(\R^d).
\end{equation} In particular, it is translation
invariant and  independent of the choice of the
convex polytope $\Delta$ on non
integer order symbols. Consequently, the map $\sigma\mapsto C(\sigma)=\cutoffsum_{\Z^d} \sigma-\cutoffint_{\R^d}\sigma$ is also
translation invariant on non integer order symbols i.e. $$C(t_{\up}^*\sigma)=
C(\sigma)\quad \forall \up\in \Z^d, \quad\forall \sigma\in CS_{\rm
  c.c}^{\notin \Z}(\R^d).$$ 
\end{thm}
\begin{rk} For a general symbol $\sigma\in CS_{\rm c.c}(\R^d)$
${\rm fp}_{N\to \infty} \sum_{N\,
  \Delta\cap \Z^d} \sigma(\un)$ depends on the choice ofthe polytope
$\Delta$. As in   \cite{MP} we choose the polytope to be a hypercube and set  $$\cutoffsum_{\Z^d}\sigma:={\rm fp}_{N\to \infty}\sum_{\vert \un\vert_{\rm sup}\leq N}
\sigma(\un)$$ which we refer to as the regularised cut-off sum of $\sigma$ on
$\Z^d$ in analogy with the similar terminology used for integrals. Note that in \cite{MP} this notation had been introduced only
for
radial functions $f(\ux)= \sigma(\vert \ux\vert_{\rm sup})$ where $\sigma$ is a
classical symbol on $\R^+$. 
\end{rk}
{\bf Proof:} Since the Hadamard finite part ${\rm fp}_{N\to \infty} \sum_{N\,
  \Delta\cap \Z^d} \sigma(\un)$ coincides with the ordinary sum $\sum_{N\in
  \Z^d} \sigma(\un)$ whenever $\sigma\in CS^{<-d}_{\rm c.c}(\R^d)$, by
Theorem \ref{thm:uniqueness}, all we need to check is the $\Z^d$-translation
invariance on non integer order symbols, i.e.  
$${\rm fp}_{N\to \infty} \sum_{N\,
  \Delta\cap \Z^d} t_{\up}^*\sigma(\un)= {\rm fp}_{N\to \infty} \sum_{N\,
  \Delta\cap \Z^d} \sigma(\un)\quad \forall \sigma\in CS_{\rm c.c}^{\notin
  \Z}(\R^d)\quad \forall \up\in \Z^d.$$
\begin{enumerate}\item Let us first show that for a symbol $\sigma$ of non
  integer order, the Hadamard finite part ${\rm fp}_{N\to
    \infty} P_{N\, \Delta}(\sigma)$ is independent of the polytope one is
  expanding. This follows from results of \cite{GSW} recalled in  Proposition
  \ref{prop:GSW} by which   the constant
$C(\sigma)$ is independent of the polytope one is expanding.  Since by Lemma \ref{lem:technicallem}, the finite parts
${\rm fp}_{R\to \infty} \left(\left(
   M^{[j]}(\partial_h)-Id\right) \tilde P_{R\, \Delta}(\sigma) (
 h)\right)_{\vert_{h=0}} $ vanish whenever the symbol has non integer order,
we infer from (\ref{eq:GSW}) that  $${\rm fp}_{R\to \infty} P_{R\, \Delta}(\sigma)= {\rm fp}_{R\to
  \infty}\tilde P_{R\, \Delta}(\sigma)+ C(\sigma)= \cutoffint_{\R^d}\sigma+
C(\sigma)$$  is independent of $\Delta$.
\item We now show translation invariance of the map $\sigma\mapsto {\rm fp}_{R\to \infty} P_{R\,
    \Delta}(\sigma)$ on non integer order symbols.
We write $\sum_{N\,
  \Delta\cap \Z^d} t_{\up}^*\sigma(\un)= \sum_{ t_{-\up}^*(N\,
  \Delta)\cap \Z^d}\sigma(\un).$ Let $N$ be an integer chosen large enough so that the
translated expanded polytope\footnote{Note that it does not coincide with the
expanded translated polytope $N\cdot t_{-\up} (\Delta)$.}
$t_{-\up} (N\Delta)$, which is  a simple regular polytope contains
$0$ in its interior. As in the case of the
original expanded polytope $N\Delta$, its  faces move out to infinity as $N\to
\infty$. We can therefore implement the same proof as in  Lemma
\ref{lem:technicallem}; using the
symbolic properties of $\sigma$ as well as the fact that it has non integer
order, we check that ${\rm fp}_{N\to \infty}
\partial^{\gamma_i}_{h_i}\left(\tilde P_{t_{-\up} (N\Delta)}(\sigma)( h)\right)_{\vert_{h=0}}=0$. By the
Khovanskii-Pukhlikov formula, we infer from there that 
$$\cutoffsum_{\Z^d}t_{-\up}^*\sigma:={\rm fp}_{N\to \infty} P_{ t_{-\up}(N\,\Delta)}(\sigma) ={\rm fp}_{N\to \infty}
\tilde P_{ t_{-\up}(N\,\Delta)}(\sigma) + C(\sigma)= \cutoffint_{\R^d}\sigma+C(\sigma)=
\cutoffsum_{\Z^d}\sigma.$$
\item Since  $\cutoffint_{\R^d}\sigma+C(\sigma)=
\cutoffsum_{\Z^d}\sigma$ and since
 both $\cutoffsum_{\Z^d}$ and $\cutoffint_{\R^d}$ are invariant under translation by
$\up\in \Z^d$ on non integer order symbols, so is the map 
$\sigma\mapsto C(\sigma)$, consequently
$$C\left(t_{\up}^*\sigma\right)= C(\sigma)\quad \forall \sigma \in CS_{\rm
  c.c}^{\notin\Z}(\R)\quad \forall \up\in \Z^d$$ 
which shows  (\ref{eq:cutoffsumnonintorder}).
\end{enumerate}
\endsquare
\vfill \eject \noindent
\section{Sums of holomorphic symbols on $\Z^d$}
Following \cite{KV} we slightly  generalise  the notion of 
   gauged symbols used in \cite{GSW} (who follow the  terminology
    introduced by V. Guillemin) in as far as we allow any non constant affine
    order whereas gauged symbols have affine order with derivative equal to
    $1$.
\begin{defn}We call local  holomorphic regularisation  of a symbol $\sigma\in CS_{\rm
    c.c.}(\R^d)$ at zero any holomorphic family\footnote{The definition of
    holomorphic families of symbols and operators  is
    recalled in the Appendix section A.4.}  $${\cal R}(\sigma): z\mapsto\sigma(z)\in CS_{\rm
    c.c.}(\R^d)$$  in a neighborhood of zero such that  $\sigma(0)=\sigma$ and $\sigma(z)$ 
  has non constant affine\footnote{This restriction is not strictly necessary
    but convenient to work with.}
order $\alpha(z)$. 
\end{defn}

\begin{ex} A Riesz perturbation ${\cal R}(\sigma)(z)(\ux):=\chi(\ux)\,  \sigma(\ux)\, \vert
  \ux\vert^z+ (1-\chi(\ux))\, \sigma(\ux)$ for some smooth cut-off function $\chi$ which vanishes in a
  neighborhood of zero and is one outside the unit ball, is a local holomorphic
  regularisation at zero. Note that $\sigma^\prime(z)(\ux)=\chi(\ux)\,
  \sigma(\ux)\, \log\vert \ux\vert\,\vert
  \ux\vert^z $ vanishes on the unit sphere. 
\end{ex}
\begin{thm}\label{thm:maintheorem} Let ${\cal R}(\sigma): z\mapsto {\cal R}(\sigma)(z):=  \sigma(z)$ be a local holomorphic
  regularisation of  $\sigma\in CS_{\rm c.c.}(\R^d)$ at zero with 
  order $\alpha(z)$ then 
the map $$z\mapsto \cutoffsum_{\Z^d}\sigma(z)(\un)$$
is meromorphic with a discrete of simple poles in  $\alpha^{-1}\left([-d, \infty[\cap
  \Z\right)$ and complex residue  at $z=0$ given by
: \begin{equation}\label{eq:discreteKV}{\rm Res}_{z=0} \cutoffsum_{\Z^d}\sigma(z)(\un)= -\frac{1}{\alpha^\prime(0)}
{\rm res}(\sigma(0)).
\end{equation}
The constant term in the Laurent series at $z=0$
$$\cutoffsum_{\Z^d}^{\cal R}\sigma(\un):=  {\rm fp}_{z=0}
\cutoffsum_{\Z^d}\sigma(z)(\un)$$  reads
\begin{equation}\label{eq:discreteregsum}
\cutoffsum_{\Z^d}^{\cal R}\sigma(\un)=\cutoffint_{\R^d}^{\cal R}\sigma(\ux)\, d\ux+C(\sigma),
\end{equation}
where we have set $\cutoffint_{\R^d}^{\cal R}\sigma(\ux)\, d\ux={\rm fp}_{z=0}
\cutoffint_{\R^d}\sigma(z)(\ux)\, d\ux.$\\
 Whenever the order of $\sigma$ has real part $<-d$ (resp. is non integer), the map $z\mapsto
 \cutoffsum_{\Z^d}\sigma(z)(\un)$ is holomorphic at $z=0$  and
 converges to 
 the ordinary sum $\sum_{\Z^d}
\sigma(\un)$ (resp. cut-off regularised sum $\cutoffsum_{\Z^d}
\sigma(\un)$)  as $z\to 0$ so that in that case
$$ \cutoffsum^{\cal R}_{\Z^d}\sigma(\un)= \sum_{\Z^d}
\sigma(\un), \quad \left({\rm resp.} \quad  \cutoffsum^{\cal R}_{\Z^d}\sigma(\un)= \cutoffsum_{\Z^d}
\sigma(\un)\right)$$
\end{thm}
\begin{rk}
Here again, $C(\sigma)$ arises as a difference of regularised integrals,
confirming a result of \cite{GSW}. \end{rk}
{\bf Proof:}  By Theorem \ref{thm:translinvregsum}, outside the set $\alpha^{-1}\left([-d, \infty[\cap
  \Z\right)$ we have:
\begin{equation}\label{eq:useful2}\cutoffsum_{\Z^d}\sigma(z)=\cutoffint_{\R^d}\sigma(z)(\ux)\, d\ux+C(\sigma(z)).\end{equation}
 On
the one hand, by results of \cite{KV} (see Theorem \ref{thm:KVsymbol} in the Appendix) we know that under
the assumptions of the theorem,  the map $z\mapsto \int_{\R^d}\sigma(z)(\ux)\,
d\ux $ is  meromorphic with a discrete set of simple poles in $\alpha^{-1}\left([-d, \infty[\cap
  \Z\right)$ and that at zero 
 \begin{equation}\label{eq:KV}{\rm Res}_{z=0} \cutoffint_{\R^d}\sigma(\ux)\, d\ux=
 -\frac{\sqrt{2\pi}^d}{\alpha^\prime(0)}{\rm res}(\sigma).\end{equation}
On the other hand, we know from \cite{GSW} that $z\mapsto C(\sigma(z))$ is
holomorphic\footnote{Their proof can easily be
  generalised to our more general setup of holomorphic families with any
  non constant affine order.}. It therefore follows from
 (\ref{eq:useful2}) that the map $z\mapsto \cutoffsum_{\Z^d}\sigma(z)(\un)$ is
 meromorphic with  a discrete of simples poles in  $\alpha^{-1}\left([-d, \infty[\cap
  \Z\right)$ and complex residue  at $z=0$ given by  $${\rm Res}_{z=0}\cutoffsum_{\Z^d}\sigma(z)(\un)=
-\frac{1}{\alpha^\prime(0)}
{\rm res}(\sigma(0)).$$
 Taking finite parts at $z=0$ in (\ref{eq:useful2}) yields
 (\ref{eq:discreteregsum}) since  $\lim_{z\to 0} C(\sigma(z))= C(\sigma).$
 \\
When the order of $\sigma$ has real part $<-d$ (resp. is non integer), the map
$z\mapsto \cutoffint_{\R^d}\sigma(z)(\ux)\, d\ux$ is holomorphic at $0$ since
$\sigma$ has vanishing residue.  Its limit at $z=0$  coincides with the ordinary
integral $\int_{\R^d}\sigma(\ux)\, d\ux$ (resp. the cut-off regularised
integral $\cutoffint_{\R^d} \sigma(\ux)\, d\ux.$). By (\ref{eq:useful2}) and
since $z\mapsto C(\sigma(z))$ is known to be holomorphic (\cite{GSW}), 
the map
$z\mapsto \cutoffsum_{\Z^d}\sigma(z)(\un)=\sum_{\Z^d}\sigma(z)(\un)$ is also holomorphic at $z=0$ and its limit
reads: 
$$ \cutoffsum_{\Z^d}^{\cal R}\sigma(\un)= \lim_{z\to 0}
\cutoffsum_{\Z^d}\sigma(z)(\un)= \int_{\R^d}\sigma(\ux)\, d\ux+ C(\sigma)=
\sum_{\Z^d}\sigma(\un),$$ 
resp. 
$$ \cutoffsum_{\Z^d}^{\cal R}\sigma(\un)= \lim_{z\to 0}
\cutoffsum_{\Z^d}\sigma(z)(\un)= \cutoffint_{\R^d}\sigma(\ux)\, d\ux+ C(\sigma)=
\sum_{\Z^d}\sigma(\un)$$where the last sum is an ordinary sum (resp. a cut-off
regularised sum). 
\endsquare
\begin{rk}Whereas regularised sums $\cutoffsum_{\Z^d}^{\cal R}$, which
  coincide with cut-off regularised sums $\cutoffsum_{\Z^d}$ on non integer
  order symbols are translation invariant on such symbols, they are not translation
  invariant on integer order symbols. For example, in dimension $d=1$, and for
non positive  negative integers  $k$ we have
$${\rm fp}_{z=0}\sum_{n\in \Z}  \vert n+p\vert^{k+z}\neq {\rm
  fp}_{z=0}\sum_{n\in \Z}  \vert n\vert^{k+z}.$$
Indeed, on the one hand,  $${\rm fp}_{z=0}\sum_{n\in \Z}  \vert n+p\vert^{k+z}=2{\rm fp}_{z=0}
\sum_{n=1}^\infty  \vert n+p\vert^{k+z}=2\, \zeta(-k,p)$$
where  $\zeta(s,p):= {\rm fp}_{z=0}
\sum_{n=1}^\infty  \vert n+p\vert^{-s+z}$ is the Hurwitz zeta function (see
e.g. \cite{C}).
On the other hand,  we
have:
   $${\rm fp}_{z=0}\sum_{n\in \Z}  \vert n\vert^{k+z}=2{\rm fp}_{z=0}
\sum_{n=1}^\infty  \vert n\vert^{k+z}=2\, \zeta(-k),$$
but $\zeta(-k)= -\frac{B_{k+1}}{k+1}\neq \zeta(-k, p)=
-\frac{B_{k+1}(p)}{k+1}$
where $B_n(x)$ stands for the $n$-th Bernouilli polynomials, and $B_n= B_n(0)$
for the $n$-th
Bernouilli constants. 
\end{rk}
\begin{cor}\label{cor:maincorollary}
Given two local holomorphic regularisations of $\sigma$, ${\cal R}$ at zero  which sends $\sigma$ to
$\sigma(z)$ of order $\alpha(z)$ and
$\tilde {\cal R}$  which sends $\sigma$ to
$\tilde\sigma(z)$ of same affine order $ \alpha(z)$  we have:
\begin{equation}\label{eq:comparison}\cutoffsum_{\Z^d}^{\cal
  R}\sigma(\un)-\cutoffsum_{\Z^d}^{\tilde{\cal R}}\sigma(\un)
=-\sqrt{2\pi}^d\,\frac{{\rm res}\left(\sigma^\prime(0)-
\tilde\sigma^\prime(0)\right)}{\alpha^\prime(0)}.
\end{equation}
In particular, 
whenever ${\cal R}(\sigma)$ and   $\tilde {\cal R}(\sigma)$  coincide with
$\sigma$  on the unit sphere, their cut-off regularised sums coincide:
$$\cutoffsum_{\Z^d}^{\cal
  R}\sigma(\un)=\cutoffsum_{\Z^d}^{\tilde{\cal R}}\sigma(\un).$$
\end{cor}
\begin{rk}A priori, as explained in the Appendix, the first derivatives $\sigma^\prime(0)$ and $\tilde \sigma^\prime(0)$ are not
  classical but log-polyhomogeneous of log-type 1 (see e.g. \cite{L},
  \cite{PS}). However, since they have same order, their difference is
  classical, so that the noncommutative residue of their difference indeed makes sense.
\end{rk}
{\bf Proof:}
  $ \tau(z):=\frac{ \sigma(z)-
\tilde\sigma(z)}{z}$ defines a local
holomorphic family of classical symbols around zero of order $\alpha(z)$ such that
$\tau(0)=\sigma^\prime(0)-
\tilde\sigma^\prime(0) $. Applying  (\ref{eq:discreteKV}) to this family $ \tau(z)$ yields
$${\rm fp}_{z=0} \cutoffsum_{\Z^d}\left(\sigma(z)(\un)-
\tilde\sigma(z)(\un)\right)=\sqrt{2\pi}^d\, {\rm Res}_{z=0} \cutoffsum_{\Z^d}\tau(z)(\un)=-\frac{{\rm res}\left(\sigma^\prime(0)-
\tilde\sigma^\prime(0)\right)}{\alpha^\prime(0)} .$$ 
If $\sigma(z)$ and $\tilde \sigma(z)$ restricted to the unit sphere are
independent of $z$, ${\rm res}\left(\sigma^\prime(0)-
\tilde\sigma^\prime(0)\right)$ vanishes and ${\rm fp}_{z=0} \cutoffsum_{\Z^d}\sigma(z)(\un)={\rm fp}_{z=0} \cutoffsum_{\Z^d}
\tilde\sigma(z)(\un)$.
\endsquare\vfill \eject \noindent
\section{Zeta functions associated with quadratic forms}Cut-off regularised
sums $\cutoffsum_{\Z^d}$ are useful  to build meromorphic
extensions of ordinary sums of holomorphic families of symbols; we recover
this way the existence of  meromorphic extensions of zeta functions associated
with quadratic forms. \\
To a  positive definite  quadratic
form $q(x_1,\cdots, x_d)$ and  a smooth cut-off function $\chi$ which vanishes in a small
neighborhood of $0$ and is identically one outside the unit euclidean ball, we
assign the classical symbol 
\begin{equation}\label{eq:quadratic}\ux\mapsto \sigma_{q, s}(\ux):=\chi(\ux)\, q( \ux)^{-s}\in
  CS_{\rm c.c}(\R^d).
\end{equation} 
\begin{thm}  \label{thm:quadraticzeta} Given any complex number 
 $s$ the  map $$z\mapsto  \sum_{\un\in \Z^d-\{0\}}\sigma_{q, s+z}= 
\sum_{\un\in \Z^d-\{0\}} q(\un)^{-(s+z)}$$ which is holomorphic on the
half plane ${\rm Re}(s+z)>d/2$ extends to a meromorphic map
$$z\mapsto  \cutoffsum_{\un\in \Z^d-\{0\}}\sigma_{q, s+z}= 
\cutoffsum_{\un\in \Z^d-\{0\}} q(\un)^{-(s+z)}$$ with
simple pole at $z=0$ given by:
$${\rm Res}_{z=0}
\cutoffsum_{\un\in \Z^d-\{0\}} q(\un)^{-(s+z)}=\delta_{2s=d}\, 
\int_{\vert\omega\vert = 1}q(\omega)^{-d/2}\, d\mu_S(\omega) $$  and constant term at $z=0$: 
\begin{equation}\label{eq:quadratic}Z_q(s):={\rm fp}_{z=0}\cutoffsum_{\un\in \Z^d-\{0\}}
q(\un)^{-(s+z)}.\end{equation}
Moreover,  
$$ Z_q(s)= \cutoffint_{\R^d} \sigma_{q,s}+ C\left(\sigma_{q,
    s}\right).$$
\end{thm}
{\bf Proof:} Up to  the pole which we compute separately, the result follows from Theorem \ref{thm:maintheorem} applied to
$\sigma_{q, s}$ and Riesz regularisation  ${\cal R}: \sigma \mapsto
\sigma(\ux)\, \vert
\ux\vert^{-z}$ combined with the fact that Riesz regularised integrals
 coincide with ordinary cut-off regularised integrals
 (see e.g. \cite{MP}, \cite{PS},  \cite{P1}):
$${\rm fp}_{z=0}\cutoffint_{\R}\sigma(\ux) \, \vert \ux\vert^{-z}\, d\ux=
\cutoffint_{\R}\sigma(\ux) \, d\ux\quad \forall\sigma
\in CS_{\rm c.c}(\R^d). $$   Now, by 
(\ref{eq:useful2}) the pole at $z=0$ is given by the pole of
$\cutoffint_{\R^d} \sigma_{q,s}(\ux)\,\vert \ux\vert^{-z}\, d\ux$. We write 
\begin{eqnarray}\label{eq:useful4}
&{}&{\rm Res}_{z=0}\cutoffint_{\R^d}\chi(\ux)\,  q(\ux)^{-s}\,\vert \ux\vert^{-z}\, d\ux\nonumber\\
&=&
{\rm Res}_{z=0}\int_{0\leq \vert\ux\vert \leq 1}\chi(\ux)\, q(\ux)^{-s}\,\vert \ux\vert^{-z}\,
d\ux +{\rm Res}_{z=0}\left({\rm
  fp}_{R\to \infty}  \int_{1\leq \vert\ux\vert \leq R}q(\ux)^{-s}\,\vert \ux\vert^{-z}\, d\ux \right)\nonumber\\
&=& {\rm Res}_{z=0}\left({\rm
  fp}_{R\to \infty}  \int_{\vert\omega\vert =1}\int_1^Rq(r\omega)^{-s}\,r^{-z+d-1}\,d\mu_S(\omega)\right) \nonumber\\
&=& 
{\rm Res}_{z=0}\left(\left({\rm
  fp}_{R\to \infty} \int_1^Rr^{-(2s+z)+d-1} \, dr\right) \, \left(
\int_{\vert\omega\vert = 1}q(\omega)^{-s}\, d\mu_S(\omega)\right)\right)\nonumber\\
&=&{\rm Res}_{z=0}\left({\rm
  fp}_{R\to \infty}\frac{ R^{-(2s+z)+d}-1}{-(s+z)+d} \right) \, \left(
\int_{\vert\omega\vert = 1}q(\omega)^{-s}\, d\mu_S(\omega)\right)\nonumber\\ &=& {\rm Res}_{z=0}\frac{1}{2s+z-d}\,  \left(
\int_{\vert\omega\vert = 1}q(\omega)^{-s}\, d\mu_S(\omega)\right)\nonumber\\
&=& \delta_{2s-d}\,  \left(
\int_{\vert\omega\vert = 1}q(\omega)^{-s}\, d\mu_S(\omega)\right).
\end{eqnarray}
As announced, there is therefore a  pole at $z=0$  only if $s=d/2$ in which
case the residue coincides with  $ \int_{\vert\omega\vert = 1}q(\omega)^{-s}\, d\mu_S(\omega)$.\endsquare
\begin{rk} For $d=2$ and $q(x,y)= ax^2+bxy+cy^2$ with $4ac-b^2>0$, 
  $Z_q(s)$ yields a meromorphic extension of   Epstein's
  $\zeta$-function 
$\sum_{(m, n)\in \Z-\{0\}} (am^2+bmn+cn^2)^{-s}$ (see e.g.\cite{CS}) which is
known to satisfy a functional equation similar to the one satisfied by the Riemann
zeta function. \\
When $a=c=1, b=0$, $Z_q(s)$ provides a meromorphic extension of the
zeta function of $\Z[i]$ given by (see e.g. \cite{C})
$$Z_4(s):= \sum_{z\in
  \Z[i]-\{0\}} \vert z\vert^{-2s}= \sum_{(m,n)\in \Z-\{0\}} m^2+n^2.$$ 
When  $a=b=c=1$, $Z_q(s)$ provides a meromorphic extension of the
zeta function of $\Z[j]$ given by (see e.g. \cite{C})
$$Z_3(s):= \sum_{z\in
  \Z[j]-\{0\}} \vert z\vert^{-2s}= \sum_{(m,n)\in \Z-\{0\}} m^2+mn +n^2.$$ 
\end{rk}
\begin{prop}\label{prop:zetaqhol} Whenever  ${\rm Re}(s)\leq 0$ 
\begin{enumerate}
\item $Z_{q}(s)=C(\ux\mapsto q(
  \ux)^{-s}),$ \item Specifically, for any non negative integer $k$
  $$Z_{q}(-k)=0.$$
\item Moreover, $Z_q$ is holomorphic at $s=-k$ for any non negative integer
  $k$ and $Z^\prime(-k)= {\partial_s}_{\vert_{ s=-k}}C(q^{-s})$
  where the derivative at $k=0$  stands for the derivative of the map
  $C(q^{-s})$ restricted to the half plane Re$(s)\leq 0$ \footnote{In contrast
    to the value at $s=-k$ which vanishes, one does not expect the derivative  to vanish in general.} . 
\end{enumerate}
\end{prop} {\bf Proof:}
\begin{enumerate}
\item When ${\rm Re}(s)\leq 0$, the map $\ux\mapsto q( \ux)^{-s}$ can be
extended by continuity  to $\ux=0$ by\footnote{ Note that this extension is not
  smooth  at $0$ so that
it does not define a symbol. It nevertheless  has the same asymptotic behaviour as
$\vert\ux\vert\to \infty$ as  $\ux\mapsto \chi(\ux)\, q( \ux)^{-s}$ which is enough
for our needs.}  $$\bar \sigma_{q,s}(\ux):= q(
\ux)^{-s}\, \forall \ux\neq 0, \quad\bar \sigma_{q,s}(0)=0.$$ In that case, there
is no need to introduce a cut-off function  $\chi$ at $0$ and we write:
$$Z_{q}(s)= {\rm fp}_{z=0,{\rm Re}(s+z)\leq 0}\cutoffsum_{\Z^d}
q(\un)^{-(s+z)}=\cutoffint_{\R^d}q(\ux)^{-s}\, d\ux+C(\ux\mapsto
q(\ux)^{-s}) $$
along the lines of the proof of the previous proposition. Here we take the
finite part at $z=0$ of the restriction $z\mapsto \cutoffsum_{\Z^d}
q(\un)^{-(s+z)}$ to the half plane Re$(z)\leq 0$. 
Using polar coordinates 
 $\ux= r\, \omega$ with $r>0$ and $ \omega$ in the unit sphere, the result then follows   from the fact that the cut-off regularised integral
vanishes if Re$(z)\leq 0$  since we have
\begin{eqnarray*}
\cutoffint_{\R^d}q(\ux)^{-2s}\, d\ux&=& {\rm
  fp}_{R\to \infty}  \int_{\vert\ux\vert \leq R}q(\ux)^{-2s}\, d\ux \\
&=& {\rm
  fp}_{R\to \infty}  \int_{\vert\omega\vert =1}\int_0^Rq(r\omega)^{-2s}\,r^{d-1}\, d\ux \\
&=&\left({\rm
  fp}_{R\to \infty} \int_0^Rr^{-2s+d-1} \, dr\right) \, \left(
\int_{\vert\omega\vert \leq 1}q(\omega)^{-2s}\, \mu_S(\omega)\right)\\
&=&\left({\rm
  fp}_{R\to \infty}\frac{ R^{-2s+d}}{-2s+d} \, dr\right) \, \left(
\int_{\vert\ux\vert \leq 1}q(\omega)^{-2s}\, \mu_S(\omega)\right)\quad
{\rm if}\quad -2s+d\neq 0\\ &=& 0\quad
{\rm if}\quad -2s+d\neq 0.
\end{eqnarray*}
where $\mu_S$ is the volume measure on the unit sphere induced by the canonical
measure on $\R^d$.\item 
When   $s=-k$, we also have $ C(\ux\mapsto
q(\ux)^k)=0$ since $C$ vanishes on polynomials so that $Z_q(-k)=0$.
\item By Theorem \ref{thm:quadraticzeta} there is no pole at $s=-k$ (the
  presence of the cut-off function $\chi$ does not affect  poles) since the
  only pole corresponds to $s=d/2$. The map $Z_q$ is therefore holomorphic at
  $s=-k$ with derivative given by the derivative of the map $C(\ux\mapsto q(
  \ux)^{-s})$ at $s=-k$. 
\end{enumerate}
\endsquare\vfill \eject  \noindent
\section{Zeta functions associated with Laplacians on tori}
  The Laplacian $-\sum_{j=1}^d\partial_j^2$  on
$\R^d$ induces  a non negative elliptic differential operator 
$\Delta_d$ of order $2$  with positive leading symbol on  $\T^d\simeq
\R^d/(2\pi\Z)^d$, the spectrum of which reads:
$${\rm Spec}(\Delta_d)=\{n_1^2+\cdots +n_d^2, \quad n_i\in \Z\}.$$
Let $\Delta_d^{^\perp}$ denote the restriction of  $\Delta_d$ to the
orthogonal of its kernel, which we recall is finite dimensional. Recall (see
e.g \cite{KV})  that the map
$z\mapsto {\rm tr}^\prime (\Delta_d^{-s}):=
{\rm tr} \left(\left(\Delta_d^\perp\right)^{-s}\right)$ which is  holomorphic
on the half-plane   Re$(s)>\frac{d}{2}$, extends to a meromorphic map 
$z\mapsto {\rm TR}^\prime\left(\Delta_d^{ -s}\right):= {\rm TR}
\left(\left(\Delta_d^\perp\right)^{-s}\right)
$ on the whole complex plane (see Theorem \ref{thm:KVop} in the Appendix), where TR stands
for Kontsevich and Vishik's  canonical trace \cite{KV}, the definition of which is
given in the Appendix (see Definition \ref{defn:KV}). The zeta
function associated with $\Delta_d$ is defined by
$$\zeta_{\Delta_d}(s):= {\rm fp}_{z=0} {\rm
  TR}^\prime\left(\Delta_d^{-(s+z)}\right).$$
Applying Theorem \ref{thm:quadraticzeta} to the quadratic form $q(\ux)= \sum_{i=1}^dx_i^2$ leads
to a description of the $\zeta$-function associated with the Laplacian as a
regularised discrete sum of powers of its eigenvalues.
\begin{prop}\label{prop:zetaDelta}\begin{equation}\label{eq:zetad}\zeta_{\Delta_d}(s)= {\rm fp}_{z=0} \cutoffsum_{\Z^d-\{0\}}
  (n_1^2+\cdots +n_d^2)^{-(s+z)}=\cutoffint_{\R^d}\sigma_{s}+
  C(\sigma_s)\end{equation}
where we have set $\sigma_s(\ux)= \chi(\ux)\, \vert \ux\vert^{-2s}$. \\
 Here
$\chi$ is any smooth cut-off function which vanishes in a  small neighborhood of
$0$ and is identically equal to $1$ outside the unit euclidean ball. 
\end{prop}
{\bf Proof:}
For ${\rm Re}(z)$ sufficiently large, we have  $${\rm tr}^\prime\left(\Delta_d^{-(s+z)}\right)=
\sum_{\Z^d-\{0\}}  (n_1^2+\cdots +n_d^2)^{-(s+z)}.$$
By theorem \ref{thm:quadraticzeta} the map  $\sum_{\Z^d-\{0\}}
(n_1^2+\cdots +n_d^2)^{-(s+z)}$   extends to a  meromorphic map $\cutoffsum_{\Z^d-\{0\}}
(n_1^2+\cdots +n_d^2)^{-(s+z)}$ which by the uniqueness of the meromorphic
extension therefore  coincides with 
 ${\rm  TR}^\prime\left(\Delta^{-(s+z)}\right)$. \\
Again by  Theorem \ref{thm:quadraticzeta} we further have
$$\zeta_{\Delta_d}(s)= \cutoffint_{\R^d}\sigma_s(\ux) \, d\ux+C(\sigma_s) $$ as anounced in  the proposition.
   \endsquare \\ \\
By  Proposition \ref{prop:zetaqhol} the  $\zeta$-function
$ \zeta_{\Delta_d}$ is  holomorphic at any non positive integers $-k$.
In particular, the zeta determinant 
$${\rm det}_\zeta(\Delta_d):= e^{-\zeta_{\Delta_d}^\prime(0)}$$ is
well-defined. The following proposition relates it to the
derivative at zero of the map  $s\mapsto C(\ux\mapsto \vert \ux\vert^{-2s})$. 
\begin{prop}\label{prop:zetadeterminant} Whenever  ${\rm Re}(s)\leq 0$ $$\zeta_{\Delta_d}(s)=C(\ux\mapsto \vert
  \ux\vert^{-2s}).$$ Specifically, for any non negative integer $k$
  $$\zeta_{\Delta_d}(-k)=0.$$
Moreover, 
$${\rm det}_\zeta(\Delta_d)= \exp\left(-{\partial_s}_{\vert_{ s=0^-}}C(\ux\mapsto \vert
  \ux\vert^{-2s})\right)$$
where the subscript ${\partial_s}_{\vert_{ s=0^-}}$ stands for the derivative of the map $s\mapsto C(\ux\mapsto \vert  \ux\vert^{-2s})$ restricted to 
 the half plane Re$(s)\leq 0$. 
\end{prop} {\bf Proof:} This follows from Proposition \ref{prop:zetaqhol}
applied to $q(\ux)= \vert \ux\vert^2$.\endsquare
\begin{ex} When $d=1$ this yields:
$$\zeta(s)= \frac{1}{2} \zeta_{\Delta_1}(-s/2)=\frac{C(x\mapsto \vert
  x\vert^{-s})}{2} \quad {\rm  for}\quad {\rm Re}(s)\leq 0.$$ In particular, for
$s=-k$ with $k$ a non negative integer, we have that
$\zeta(-k)= \frac{1}{2} \zeta_{\Delta_1}(-k/2)=\frac{C(x\mapsto \vert
  x\vert^k)}{2}$ vanishes for even $k$. 
\end{ex} 
When ${\rm Re}(s)>d/2$, the sum $\sum_{\Z^d} \vert\un\vert^{-2s}\, \log\vert
\un\vert$ converges and  we have
$$\zeta_{\Delta_d}^\prime(s)=-2\, \sum_{\Z^d-\{0\}} \vert\un\vert^{-2s}\, \log\vert
  \un\vert\quad{\rm if}\quad{\rm Re}(s)\leq d/2.$$
Just as one can extend cut-off regularised integrals to log-polyhomogeneous
symbols,   cut-off regularised sums can be extended to
log-polyhomogeneous symbols so that one can define $\cutoffsum_{\Z^d-\{0\}} \vert\un\vert^{-2s}\, \log\vert
  \un\vert$ for general complex numbers $s$.
But {\it in general}, 
$$\zeta_{\Delta_d}^\prime(s)\neq-2\, \cutoffsum_{\Z^d} \vert\un\vert^{-2s}\, \log\vert
  \un\vert\quad{\rm when}\quad{\rm Re}(s)\leq d/2.
$$
\begin{rk}A similar observation holds for the $\zeta$ function itself for in general
$\zeta_{\Delta_d}(s)\neq \cutoffsum_{\Z^d} \vert\un\vert^{-2s},$ as
 for example when $s=-k$ is a non positive integer. In that case, the right
 hand side vanishes whereas the left hand side does not. 
\end{rk}
Consequently, even though the map $z\mapsto C(\sigma(z))$ is holomorphic at
$z=0$ for any local holomorphic regularisation ${\cal R}(\sigma): z\mapsto
\sigma(z)$ of a symbol $\sigma\in CS_{\rm
  c.c}(\R^n)$ so that $\lim_{z\to 0}C(\sigma(z))= C(\sigma)$,  one  expects
that  for $\sigma$ of order $a$, in {\it general}
$${\partial_z}_{\vert_{z=0}} C(\sigma(z))\neq C(\sigma^\prime(0))\quad{\rm
  when}\quad{\rm Re}(a)\geq- d.$$
\begin{rk}Comparing the  zeta-regularised determinant
  $\zeta_{\Delta_d}^\prime(0)$ with the ``cut-off regularised'' determinant
$\exp\left(-2\, \sum_{\Z^d-\{0\}} \, \log\vert
  \un\vert\right)$ is reminiscent of issues addressed in \cite{FG} in the
context of 
  Szeg\"o operators. 
\end{rk}
\begin{ex}It follows from results of  \cite{MP} that the difference
$\zeta_{\Delta_1}^\prime (-k/2)+2\ \cutoffsum_{\Z} \vert n\vert^{k}\, \log\vert
  n\vert$ is a rational number for any positive integer $k$.
\end{ex}

\section*{Appendix : Prerequisites on regularised integrals of symbols }
This appendix recalls  known results  on regularised integrals of
classical symbols used in the main bulk  of the paper. 
\subsection*{A.1 Classical  pseudodifferential symbols with constant coefficients}
 We only give a few definitions and refer the reader to
\cite{Sh,T,Tr} for further details on classical pseudodifferential
symbols.\\  For any complex number $a$, let us denote by ${\cal S}^a_{\rm
  c.c}(\R^n)$ 
 the set of smooth functions on $\R^n$ called symbols with constant
 coefficients, such that for any multiindex
$\beta\in \N^n$ there is a constant $C(\beta)$ satisfying the following requirement:
$$\vert\partial_\ux^\beta \sigma(\ux)\vert\leq C(\beta) \vert (1+\vert
\ux\vert)^{{\rm Re}(a)-\vert \beta\vert}$$ where Re$(a)$ stands for the real
part of $a$, $\vert
\ux\vert$ for the euclidean norm of $\ux$. We single out the subset  $CS^a_{\rm c.c}(\R^n)\subset {\cal
  S}_{\rm c.c}^a(\R^n) $ of   symbols $\sigma$, called classical symbols of order $a$  with constant
coefficients,   such that 
\begin{equation}\label{eq:asymptsymb}
\sigma(\ux)=
\sum_{j=0}^{N-1} \chi(\ux)\, \sigma_{a-j}( \ux) +\sigma_{(N)}( \ux)
\end{equation}
where $\sigma_{(N)}\in {\cal S}_{\rm c.c}^{a-N}(\R^n)$ and
 where $\chi$ is a smooth cut-off function which vanishes in a
  small
   ball of $\R^n$ centered at $0$ and which is constant equal to $1$ outside
   the unit ball. Here $\sigma_{a-j,}, j\in \N_0$ are positively homogeneous of degree
   $a-j$. \\The ordinary product of functions sends $CS^a_{\rm c.c}(\R^n )\times
CS^b_{\rm c.c}(\R^n)$ to $CS^{a+b}_{\rm c.c}(\R^n)$ provided $b-a\in \Z$; let   
\begin{equation}\label{eq:CSRn} CS_{\rm c.c}(\R^n)= \langle
\bigcup_{a\in \C}CS^a_{\rm c.c}(\R^n)\rangle 
\end{equation} denote the algebra generated by   all classical
symbols with constant coefficients  on $\R^n$.  Let  $$CS_{\rm c.c}^{-\infty}(\R^n)=  \bigcap_{a\in
  \C}CS^a_{\rm c.c}(\R^n)$$ be the algebra of smoothing symbols. We write
$\sigma\sim \sigma^\prime$ for two symbols $\sigma, \sigma^\prime$  which differ by a smoothing symbol.\\
 We also
denote by $CS_{\rm
  c.c}^{<p}(\R^n):= \bigcup_{{\rm Re}(a)<p} CS_{\rm c.c}^a(\R^n)$, the set of classical
symbols of order with real part $<p$ and by   \begin{equation}\label{eq:CSnoninteger}CS_{\rm
  c.c}^{\notin \Z}(\R^n):= \bigcup_{a\in\C- \Z} CS_{\rm c.c}^a(\R^n)
\end{equation}
the set of non integer order symbols.
\subsection*{A. 2 The noncommutative residue and cut-off regularised  integrals on
  symbols }
We first
recall the definition of the noncommutative residue of a classical symbol
\cite{G2}, \cite{W1,W2}.
\begin{defn}  The noncommutative residue  is a linear form on  $CS_{\rm c.c}(\R^n)$ defined by
$${\rm res}(\sigma):=\frac{1}{\sqrt{2\pi}^n}\,  \int_{S^{n-1}}\sigma_{-n}(\ux)\,
d\mu_S(\ux)$$
where  $$d \mu_S( \ux):= \sum_{j=1}^n (-1)^{j-1} \,
\ux_j\,d\ux_1\wedge \cdots \wedge d \hat \ux_j\wedge\cdots \wedge d\ux_n$$
denotes     the  volume measure on the unit sphere $S^{n-1}$ induced by the
canonical measure on $\R^n$.
\end{defn}Let us now recall the construction of a useful linear  extension of the ordinary integral.\\
  For any $R>0$,  $B(0, R)$
denotes the ball of radius $R$ centered at $0$ in $\R^n$.
We recall that given a symbol $\sigma\in CS^a_{\rm c.c}(\R^n)$,  the map
$R\mapsto \int_{B(0, R)} \sigma (\ux)\,d\,  \ux$ has an asymptotic
expansion as $R\to\infty$ of  the form (with  the notations of
(\ref{eq:asymptsymb})):\begin{equation}\label{eq:Rasymptotics}\int_{B(0, R)}\sigma(\ux) \, d\, \ux\\
\sim_{ R\to \infty}\alpha_0(\sigma)+
\sum_{j=0,a-j+n\neq 0}^\infty  \sigma_{a-j }\,  R^{a-j+n}+
{\rm res}(\sigma) \cdot  \log  R.
\end{equation} 
\begin{defn}\label{defn:cutoffint} Given $\sigma\in CS^a_{\rm c.c}(\R^n)$ with
  $a\in \C$, we call the constant term
  $\alpha_0(\sigma)$ the cut-off regularised integral of $\sigma$:
$$\cutoffint_{\R^n}  \sigma (\ux) \, d\,\ux:= {\rm fp}_{R\to \infty} \int_{B(0, R)}\sigma(\ux) \, d\,
  \ux.$$
\end{defn}
The cut-off regularised integral $\cutoffint_{\R^n}$, which reads
\begin{eqnarray}\label{eq:constanttermclassical}
\cutoffint_{\R^n}  \sigma (\ux) \, d\,\ux& =&  \int_{\R^n}  \sigma_{(N)} (\ux)\, d\,\ux+\sum_{j=0}^{N-1}  \int_{B(0, 1)}
 \chi(\ux)\,   \sigma_{a-j}(\ux) \, d\, \ux\nonumber\\
&  -& \sum_{j=0, a-j+n\neq 0}^{N-1}  \frac{ 1}{a-j+n}
 \int_{S^{n-1}} \sigma_{a-j } (\omega)\,d  \mu_S(\omega).
\end{eqnarray}  defines a linear form on
$CS_{\rm
  c.c}(\R^n)$ which extends the ordinary integral in the following sense; if
$\sigma$ has complex order with real part smaller than $-n$ then $\int_{B(0, R)} \sigma(\ux)\,
d\, \ux$ converges as $R\to \infty$ and $$\cutoffint_{\R^n} \sigma(\ux)\,
d\,\ux=\int_{ \R^n}  \sigma (\ux)\, d\,\ux.$$
\subsection*{A. 3 The noncommutative residue and canonical trace on operators}
Let $U$ be a connected open subset of $\R^n$. 
 \\ For any complex number $a$, let ${\cal S}_{\rm cpt}^a(U)$ denote the set of smooth
 functions on $U\times \R^n$ called symbols with compact support in $U$, such that for any multiindices
$\beta, \gamma\in \N^n$,  there is a constant $C(\beta, \gamma)$ satisfying the following requirement:
$$\vert\partial_\xi^\beta\partial_x^\gamma \sigma(x,\xi)\vert\leq C(\beta, \gamma) \vert (1+\vert
\xi\vert)^{{\rm Re}(a)-\vert \beta\vert}$$ where Re$(a)$ stands for the real
part of $a$, $\vert
\xi\vert$ for the euclidean norm of $\xi$. We single out the subset
$CS^a_{\rm cpt}(U)\subset {\cal
  S}_{\rm cpt}^a(U) $ of   symbols $\sigma$, called classical symbols of
order $a$  with compact support in $U$,   such that 
\begin{equation}\label{eq:asymptsymbbis}
\sigma(x,\xi)=
\sum_{j=0}^{N-1} \chi(\xi)\, \sigma_{a-j}(x, \xi) +\sigma_{(N)}(x, \xi)
\end{equation}
where $\sigma_{(N)}\in {\cal S}_{\rm cpt}^{a-N}(U)$ and
 where $\chi$ is a smooth cut-off function which vanishes in a
  small
   ball of $\R^n$ centered at $0$ and which is constant equal to $1$ outside
   the unit ball. Here $\sigma_{a-j}(x, \cdot), j\in \N_0$ are positively homogeneous of degree
   $a-j$.
\\
  Let  $$CS_{\rm cpt}^{-\infty}(U)=  \bigcap_{a\in
  \C}CS^a_{\rm cpt}(U)$$ be the set of smoothing symbols with compact
support in $U$; we write $\sigma\sim \tau$ for two symbols that differ by a
smoothing symbol. \\
 The star product
\begin{equation}\label{eq:starproduct}\sigma\star \tau\sim\sum_{\alpha} \frac{(-i)^{\vert \alpha\vert}}{\alpha!} \partial_\xi^\alpha \sigma\,
\partial_x^\alpha \tau
\end{equation}
  of  symbols  $\sigma\in CS_{\rm cpt}^a(U)$ and $\tau\in CS_{\rm cpt}^b(U)$
  lies in  $CS_{\rm cpt}^{a+b}(U)$ provided $a-b\in \Z$.  \\
 Let  $$ CS_{\rm cpt}(U)= \langle
\bigcup_{a\in \C}CS^a_{\rm cpt}(U)\rangle $$ denote the algebra  generated by  all classical
symbols   with compact support in  $U$. We
denote by $CS_{\rm
  cpt}^{<p}(U):= \bigcup_{{\rm Re}(a)<p} CS_{\rm cpt}^a(U)$, the set of classical
symbols  of order with real part $<p$  with compact support in $U$, by  $CS_{\rm
  cpt}^{ \Z}(U):= \bigcup_{a\in\Z} CS_{\rm cpt}^a(U)$ the algebra  of  integer
order symbols, and  by   $CS_{\rm
  cpt}^{\notin \Z}(U):= \bigcup_{a\in\C- \Z} CS_{\rm cpt}^a(U)$
the set of non integer order symbols with compact support in $U$.\\
Both the noncommutative residue and the cut-off regularised integral extend to $CS_{\rm cpt}(U)$.
\begin{defn}
\begin{enumerate}
\item  The noncommutative residue of a symbol $\sigma \in CS_{\rm cpt}(U)$
  is defined by
$${\rm res}(\sigma):=\frac{1}{(2\pi)^n}\, \int_Udx \int_{S^{n-1}}\sigma_{-n}(x,\xi)\,
\mu_S(\xi)=\frac{1}{\sqrt{2\pi}^n}\,  \int_{U}{\rm res}_x(\sigma)\,d\, x  $$
where
${\rm res}_x(\sigma):=\frac{1}{\sqrt{2\pi}^n} \int_{ S^{n-1}}\sigma_{-n}(x, \xi)\, d\mu_S(\xi)$
 is the residue density at point $x$ and where  as before  $$ d\mu_S( \xi):= \sum_{j=1}^n (-1)^j \,
\xi_j\,d\xi_1\wedge \cdots \wedge d \hat \xi_j\wedge\cdots \wedge d\xi_n$$
denotes     the  volume measure on $S^{n-1}$ induced by the
canonical measure on $\R^n$.
\item  For any   $\sigma \in CS_{\rm cpt}(U)$ the cut-off regularised integral  of $\sigma$
  is defined by
$$\cutoffint_{T^*U}\sigma:= \int_Udx \cutoffint_{T_x^*U}\,
\sigma(x,\xi)\,d\xi.  $$
\end{enumerate}
\end{defn} 
Let $M$ be an $n$-dimensional closed connected Riemannian manifold (as before $n>1$).  For $a\in
\C$, let $\Cl^{a}(M)$ denote the
linear space of classical pseudodifferential operators of order
$a$, i.e. linear maps acting on smooth functions $\Ci(M)$, which using a
partition of unity adapted to an atlas on $M$ can be written as a finite sum of operators
$$A= {\rm Op}(\sigma(A))+ R$$
where $R$ is a linear operator with smooth kernel and 
$\sigma(A)\in CS^a_{\rm cpt}(U)$ for some open subset $U\subset
\R^n$. Here we have set
$${\rm Op}(\sigma)(u):= \int_{\R^n} e^{i\langle x-y, \xi\rangle} \sigma(x,
\xi)\, u(y)\, dy\,
d\xi$$
where $\langle\cdot, \cdot\rangle$ stands for the canonical scalar product in
$\R^n$. \\The star product (\ref{eq:starproduct}) on classical symbols with
compact support induces the operator product on (properly supported) classical pseudodifferential operators since
${\rm Op}(\sigma\star \tau)={\rm Op}(\sigma) \, {\rm Op}( \tau)$. 
It follows that the product $AB$ of two classical pseudodifferential 
operators $A\in \Cl^a(M)$, $B\in \Cl^b(M)$ lies in $
\Cl^{a+b}(M)$ provided $a-b\in \Z$. Let us
denote by $\Cl(M)=\langle \bigcup_{a\in \C}\Cl^{a}(M)\rangle$ the algebra
generated by all classical pseudodifferential operators acting on $\Ci(M)$. 
\\
  Given a finite rank vector bundle $E$ over $M$ we set
$ \Cl^a(M, E):= \Cl^a(M)\otimes {\rm
  End}(E)$, $ \Cl(M, E):= \Cl(M)\otimes {\rm
  End}(E)$.\\
The sets $  \Cl^{\in
  \Z}(M, E)$ and $ \Cl^{\notin
  \Z}(M, E)$ are defined
similarly using
trivialisations of $E$   from the sets $CS^\Z_{\rm cpt}(U)$ and
$\Cl^{\notin \Z}_{\rm cpt}(U)$.\\
Using a partition of unity, one can patch up the noncommutative residue,
resp. 
the cut-off regularised integral
of symbols with compact support to a noncommutative residue on all classical
pseudodifferential operators \cite{G1}, \cite{W1, W2}, resp. a  canonical trace on non integer order
classical pseudodifferential operators \cite{KV}.
\begin{defn}\label{defn:KV} 
\begin{enumerate}
\item 
The noncommutative residue  is defined on   $ \Cl(M,E)$ 
by
$${\rm res}(A):=\frac{1}{(2\pi)^n}\, \int_M dx \int_{S_x^*M}{\rm tr}_x\left(\sigma(A)\right)_{-n}(x,\xi)\,
\mu_S(\xi)=\frac{1}{\sqrt{2\pi}^n}\,  \int_{M}{\rm res}_x(A)\,d\, x  $$
where
${\rm res}_x(A):=\frac{1}{\sqrt{2\pi}^n} \int_{S_x^*M}{\rm tr}_x\left(\sigma(A)\right)_{-n}(x, \xi)\, \mu_S(\xi)$
is  the residue density at point $x$.
\item  The canonical trace  is  defined on $\Cl^{\notin \Z}(M, E)$ by
$${\rm TR}(A):=\frac{1}{(2\pi)^n}\, \int_Mdx \cutoffint_{T_x^*M}{\rm tr}_x\left(\sigma(A)(x, \xi)\right)\,
d\,\xi=\frac{1}{\sqrt{2\pi}^n}\,  \int_{M}{\rm TR}_x(A)\,d\, x  $$
where
${\rm TR}_x(A):=\frac{1}{\sqrt{2\pi}^n} \cutoffint_{T_x^*M }{\rm tr}_x\left(\sigma(A)(x, \xi)\right)\, d\xi$
 is the canonical trace density at point $x$.
\end{enumerate}
\end{defn}
\subsection*{A. 4  Holomorphic families of symbols and operators}
The notion of holomorphic family of classical pseudodifferential operators
used by Kontsevich and Vishik in \cite{KV} 
generalises the notion of complex power $A^z$ of an elliptic operator developped by Seeley \cite{Se}, the
derivatives of which lead to logarithms.  
\begin{defn}
 Let $\Omega$ be a domain of $\C$ and $U$ an open subset of $\R^d$. A family
 $(\sigma(z))_{z\in \Omega}\subset CS(U)$ is holomorphic when \\
(i) the order $\alpha(z)$ of $\sigma(z)$ is holomorphic  on $\Omega$.\\
(ii) The map $z\to\sigma(z) $ is holomorphic on $\Omega$ and
 $\forall k \geq 0, \partial_z^k\sigma(z)\in S^{\alpha(z)+\e}(U)$ for all $\e>0$.\\
(iii) For any integer $j\geq 0,$ the (positively) homogeneous component
$\sigma_{\alpha(z)-j}(z) $ of degree $\alpha(z)-j$ of the symbol  is holomorphic on $\Omega.$
\end{defn}
The derivative  of  a holomorphic family $\sigma(z)$ of classical
  symbols yields a holomorphic family of symbols, the  asymptotic expansions of
  which a priori involve a logarithmic term.
\begin{lem} The derivative of  a holomorphic family $\sigma(z)$ of classical
  symbols of order $\alpha(z)$ defines a
holomorphic family of symbols $\sigma^\prime(z)$  of order $\alpha(z)$ with
asymptotic expansion:
\begin{equation}\label{eq:sigmaprime}\sigma^\prime(z)( x,\xi)\sim\sum_{j=0}^\infty  \chi(\xi)\left( \log \vert
  \xi\vert\, \sigma_{\alpha(z)-j, 1}^\prime(z)( x,\xi)
+ \sigma^\prime_{\alpha(z)-j, 0}(z)( x,\xi)\right)\quad\forall (x, \xi)\in T_x^*U
\end{equation} for some smooth cut-off function $\chi$ around the origin which is
identically equal to $1$ outside the open unit ball and positively
homogeneous symbols
\begin{equation}\label{eq:sigmaprimej}\sigma_{\alpha(z)-j, 0}^\prime(z)(x, \xi)=
\vert \xi\vert^{\alpha(z)-j}\,  \partial_z
\left(\sigma_{\alpha(z)-j}(z)(\frac{\xi}{\vert
  \xi \vert})\right), \quad 
\sigma^\prime_{{\alpha(z)-j}, 1}(z)=\alpha'(z)\, 
\sigma_{\alpha(z)-j}(z) 
\end{equation}
of degree $\alpha(z)-j$. 
\end{lem}
\begin{rk} If $\sigma(z)$ is independent of $z$ then  $\sigma^\prime_{\alpha(z)-j} $  restricted to the unit
  sphere vanishes. 
 \end{rk} {\bf Proof:} We write 
$$\sigma(z)(x, \xi)\sim\sum_{j=0}^\infty  \chi(\xi) \, \sigma_{\alpha(z)-j}(z)(
x,\xi).$$ Using  the positive homogeneity of the components
$\sigma_{\alpha(z)-j}$ we have:
\begin{eqnarray*}
&{}&\partial_z \left(\sigma_{\alpha(z)-j}(z)(x, \xi)\right)\\
&=& \partial_z \left(\vert \xi
\vert^{\alpha(z)-j}\sigma_{\alpha(z)-j}(z)(x,\frac{\xi}{\vert \xi \vert}) \right)\\
&=& \left(\alpha'(z)\vert \xi \vert^{\alpha(z)-j}\sigma_{\alpha(z)-j}(z)(x,\frac{\xi}{\vert \xi
\vert}) \right)\log \vert \xi \vert
 +\vert \xi
\vert^{\alpha(z)-j} \partial_z
\left(\sigma_{\alpha(z)-j}(z)(x,\frac{\xi}{\vert
  \xi \vert})\right)\\
&=& \left(\alpha'(z)\sigma_{\alpha(z)-j}(z)(x,\xi) \right)\log \vert \xi \vert
 +\vert \xi
\vert^{\alpha(z)-j} \partial_z
\left(\sigma_{\alpha(z)-j}(z)(x,\frac{\xi}{\vert
  \xi \vert})\right)
\end{eqnarray*}
which shows that $\partial_z \left(\sigma_{\alpha(z)-j}(z)(x, \xi)\right)$ has
order $\alpha(z)-j$. 
Thus 
$$\partial_z\left( \sigma_N(z)(x, \xi)\right)= \sigma^\prime(z) ( x,\xi)- \sum_{j<N} \chi(\xi)\,
\partial_z\left(\sigma_{\alpha(z)-j}(z)(  x,\xi)\right)$$
lies in ${\cal S}^{\alpha(z)-N+\e}(U)$ for any $\e>0$ so that
$\sigma^\prime(z)$ is a symbol of order $\alpha(z)$ with asymptotic expansion:
\begin{equation}\sigma^\prime(z)(x, \xi)\sim
  \sum_{j=0}^\infty \chi(\xi)\, \sigma_{\alpha(z)-j}^\prime(z) \quad\forall ( x,\xi)\in
  T_X^*U\end{equation}
where
$$\sigma^\prime_{\alpha(z)-j}(z)( x,\xi):=\log \vert \xi\vert
 \sigma^\prime_{\alpha(z)-j,1}(z)(x, \xi)+  \sigma^\prime_{\alpha(z)-j,0}(z)(x, \xi)$$
for some positively homogeneous symbols $$\sigma^\prime_{\alpha(z)-j,0}(z)(x,
\xi):=\vert \xi\vert^{\alpha(z)-j}\,  \partial_z
\left(\sigma_{\alpha(z)-j}(z)(x,\frac{\xi}{\vert
  \xi \vert})\right)$$
and $$\sigma^\prime_{\alpha(z)-j,1}(z)(x, \xi):= \alpha^\prime(z)\,\sigma_{\alpha(z)-j}(z)(x, \xi) $$ of degree $\alpha(z)-j$. \\
 On the other hand, differentiating the asymptotic expansion $\sigma(z)(x, \xi)\sim \sum_{j=0}
\chi(\xi)\, \sigma_{\alpha(z)-j} (z)( x,\xi) $  w.r.
to $z$ yields
 $$
\sigma^\prime(z)(x, \ux)\sim \sum_{j=0} \chi(\xi)\,
\partial_z\left(\sigma_{\alpha(z)-j} (z)( x,\xi)\right). $$
Hence, 
$$\partial_z\left(\sigma_{\alpha(z)-j}(z) (
  x,\xi)\right)=\sigma^\prime_{\alpha(z)-j}(z)(x, \xi)=\vert \ux\vert^{\alpha(z)-j}\,  \partial_z
\left(\sigma_{\alpha(z)-j}(z)(x,\frac{\xi}{\vert
  \xi \vert})\right)+   \alpha^\prime(z)\,\sigma_{\alpha(z)-j}( x,\xi)\, \log
\vert \xi\vert $$ as announced.
\endsquare\\ 
The notion of holomorphic family  extends to operators as follows.
\begin{defn} A family $(A(z))_{z \in \Omega}\in \Cl(M, E)$ is holomorphic if
  in any local trivialisation   we can write $A(z)$ in the form
 $A(z)={\rm Op}(\sigma(A(z)))+R(z)$, for some holomorphic family of symbols
 $\left(\sigma(A(z))\right)_{z\in \Omega}$ and some holomorphic family $(R(z))_{z \in
   \Omega}$ of smoothing operators i.e. given by a holomorphic family of
 smooth Schwartz kernels.
\end{defn}
\subsection*{A. 5 Defect formulae for regularised integrals and traces}
The noncommutative residue deserves its name since it is proportional to a
complex residue as shows the following theorem. It also gives a ``defect
formula'' which compares the finite part at the poles $z_j$  of the meromorphic expansion given by
$z\mapsto \cutoffint_{\R^d}\sigma(z)$ with $ \cutoffint_{\R^d}\sigma(z_j)$.
 \begin{thm}\label{thm:KVsymbol} 
 Let $U$ be an open subset of $\R^n$ and let $z\mapsto \sigma(z)\in
 CS^{\alpha(z)}(U)$
be a holomorphic family of classical pseudo-differential symbols   parametrised by a domain
$\Omega\subset \C$ with non constant affine order $\alpha(z)$. Then
the map
$z\mapsto \cutoffint_{T_x^*U}\sigma(z)$
is meromorphic with   poles  of order $1$ at points $z_j\in \Omega\cap
\alpha^{-1}\left([-n, +\infty[\, \cap \, \Z\right)$ and
\begin{enumerate}
\item\cite{KV}
$${\rm Res}_{z=z_j} \cutoffint_{T_x^*U}\sigma(z)(x, \xi)\,
d\xi=-\frac{ \sqrt{2\pi}^n}{\alpha^\prime(z_j)}{\rm res}(\sigma(z_j)).$$
\item Moreover \cite{PS}, 
$${\rm fp}_{z=z_j} \cutoffint_{T_X^*U}\sigma(z)(x, \xi)\,
d\xi=\cutoffint_{T_X^*U}\sigma(z_j)-\frac{\sqrt{2\pi}^n}{\alpha^\prime(z_j)}{\rm
  res}_{x}(\sigma^\prime(z_j)).$$
\end{enumerate}
\end{thm}
\begin{rk}
Here the noncommutative residue has been extended to the a priori
log-polyhomogeneous symbol $\sigma^\prime(z_j)$ in a  straightforward manner by
the same formula as for classical symbols.
\end{rk}\begin{defn}
Given  a symbol $\sigma\in CS(U)$ we call  a local 
holomorphic family 
$$\sigma(z)( x, \xi)\sim  \sigma ( x,\xi) \, \vert
\xi\vert^{-z } $$ such that $\sigma(0)= \sigma$  a {\rm  Riesz regularisation of}
$\sigma$. 
\end{defn} 
\begin{ex}$\sigma(z)( x,\xi)= (1-\chi(\xi)) \sigma(x, \xi)+ \chi (\xi)\, \sigma ( x,\xi)\, \vert
\xi\vert^{-z }$ is a Riesz regularisation of the symbol $\sigma$, which depends on the
choice of cut-off function $\chi$ around zero. However, as we shall see later
on, this dependence does not affect the corresponding regularised integrals. 
 \end{ex}
Specializing Theorem \ref{thm:KVsymbol} to Riesz regularisations, we have:
\begin{cor}\label{cor:Riesz} 
\begin{enumerate}
\item Let $\sigma(z)$ be a Riesz
  regularisation of a symbol $\sigma \in CS(U)$. The   map 
$$z\mapsto\cutoffint_{\R^d}\sigma(z)(x, \xi) \, d\xi$$
is meromorphic with simple  poles at $z=0$.
\item The finite part at $z=0$ and  the ``cut-off'' finite part
coincide:
\begin{equation}\label{eq:cutoffRieszfinitepart} \cutoffint_{T_x^*U}
  \sigma(x, \xi)\, d\xi= {\rm fp}_{z=0}
\cutoffint_{T_x^*U} \sigma(z)(x, \xi) \,d \, \xi. 
\end{equation}
In particular, the finite part at $z=0$ is independent of the choice of
cut-off function involved in the definition of the Riesz regularisation of
$\sigma$. 
\end{enumerate}
\end{cor}
Theorem \ref{thm:KVsymbol}   extends to 
classical pseudodifferential operators.
 \begin{thm}\label{thm:KVop} 
 Let $z\mapsto A(z)\in
 \Cl^{\alpha(z)}(M, E)$
be a holomorphic family of classical pseudo-differential operators  with non constant affine order $\alpha(z)$. Then
the map
$z\mapsto {\rm TR}(A(z))$
is meromorphic with   poles  of order $1$ at points $z_j\in \Omega\cap
\alpha^{-1}\left([-n, +\infty[\, \cap \, \Z\right)$ and
\begin{enumerate}
\item\cite{KV}
$${\rm Res}_{z=z_j}{\rm TR}(A(z)) \,
dz=-\frac{1}{\alpha^\prime(z_j)}{\rm res}(A(z_j)).$$
\item Moreover \cite{PS}, 
$${\rm fp}_{z=z_j}{\rm TR}(A(z)) =\int_Mdx\left(\cutoffint_{\R^d}\sigma(A(z_j))-\frac{1}{\alpha^\prime(z_j)}{\rm
  res}(\sigma(A^\prime(z_j)))\right).$$
\end{enumerate}
\end{thm}
\vfill \eject \noindent
\bibliographystyle{plain}

\end{document}